\documentclass[11pt]{amsart}
\usepackage{amssymb, amsmath}
\usepackage{mathtools}
\usepackage{graphicx}
\usepackage{color}
\usepackage{float}
\usepackage{etex}
\usepackage[all]{xy}
\usepackage[dvips]{epsfig}
\usepackage{graphics}
\usepackage{amsbsy, amsmath,amsfonts,amssymb,amsthm,amscd,amssymb,latexsym}

\usepackage{ascmac}

\newcommand{\U}{\mathcal{U}}
\newcommand{\V}{\mathcal{V}}

\def\b0{\mbox{\boldmath $0$}}

\def\I{{\rm I}}

\def\V{{\rm V}}

\newtheorem{thm}{\bf Theorem}[section]
\newtheorem{cor}[thm]{\bf Corollary}
\newtheorem{lem}[thm]{\bf Lemma}
\newtheorem{prop}[thm]{\bf Proposition}
\newtheorem{definition}[thm]{\bf Definition}

\newtheorem{rem}[thm]{\bf Remark}
\newtheorem{exam}[thm]{\bf Example}

\begin{document}
\title[
Existence of periodic orbits for $\mathsf{PSVF}$ via Conley theory.
]
{Existence of periodic orbits for piecewise-smooth vector fields with sliding region via Conley theory.\\
}
\author[A.~Romero]{Angie T. S. Romero}
\address[A.~Romero]
{Instituto de Matem\'atica e Estat\'istica\\
Universidade Federal de Goi\'as\\
Goi\^ania-GO\\
Brazil.
}
\email{angieromero@discente.ufg.br}

\author[E.~Vieira]{Ewerton R. Vieira}
\address[E.~Vieira] 
{The Center for Discrete Mathematics and Theoretical Computer Science\\
	Rutgers University\\
	Piscataway, New Jersey\\
	USA and Instituto de Matem\'atica e Estat\'istica\\
	Universidade Federal de Goi\'as\\
	Goi\^ania-GO\\
	Brazil,
}
\email{ewerton.v@rutgers.edu, ewerton@ufg.br}

%
%
%
%

%
%
\begin{abstract}{
The Conley theory has a tool to guarantee the existence of periodic trajectories in isolating neighborhoods of semi-dynamical systems.  We prove that the positive trajectories generated by a piecewise-smooth vector field $Z=(X, Y)$ defined in a closed manifold of three dimensions without the scape region produces a semi-dynamical system. Thus, we have built a semiflow that allows us to apply the classical Conley theory. Furthermore, we use it to guarantee the existence of periodic orbits in this class of piecewise-smooth vector fields.  
}
\end{abstract}

\maketitle

\setlength{\baselineskip}{14pt}

\section{Introduction}

A dynamical system describes the evolution of a phenomenon over time, and these may be considered in discrete or continuous time. Traditionally, the solution of a differential equation system model is a continuous dynamical system.   Usually, it is not easy to explicitly obtain solutions of a dynamical system; the existence of periodic orbits and limit cycles (isolated periodic orbits). At the beginning of the 19th century, Henri Poincaré, while investigating the movement of planets, established the modern dynamic systems; his work combined topology and geometry, thus conducting a qualitative study of dynamic systems. At the end of the 19th century, David Hilbert, at the Second International Congress of Mathematicians held in Paris, proposed a list of 23 relevant problems. One famous unsolved problem is the 16th on this list, which refers to finding the maximum number of limit cycles of a field of polynomial vectors of degree greater than or equal to $2$.
\newline

In dynamical systems, a modern research theme refers to the piecewise smooth vector fields ($\mathsf{PSVF}$ for short).  The $\mathsf{PSVF}$ are systems that are not completely differentiable but are differentiable by parts, where a vector field is suddenly interrupted and changed by another distinct vector field. These systems are widely used for modeling some problems associated with control theory, economics, and biology, see \cite{bernardo}. We are interested in guaranteeing the existence of periodic orbits in $\mathsf{PSVF}$. So, the main objective of this paper is to use the Conley theory to obtain periodic orbits in $\mathsf{PSVF}$.
\newline

The Conley theory has a significant quantity of applications for the study of the dynamical and semi-dynamical systems. Charles Conley introduced a new topological index, called Conley Index, as a generalization of the Morse Index. The Conley Index guarantees the existence of invariant sets within a particular compact; one of these invariant sets may be a periodic orbit, which is one of our main focuses. In \cite{Conley}, Charles Conley began developing this theory for two-sided flows on compact or locally compact spaces and this was continued by Dietmar Salamon, see \cite{Salamon}, and was extended to semiflows by Rybakowski, see \cite{Rybakowski}. In the paper \cite{Mccord}, the authors present a result towards finding periodic orbits in semi-dynamical systems; it is the main tool for this paper. A modern application of finding periodic orbit in neuroscience is in the thesis \cite{Abel}, where the author uses numerical techniques to obtain periodic orbits in dynamics given by Competitive Threshold-Linear Networks and Wilson-Cowan networks.
\newline

Since the periodic orbits do not correspond to local objects, our study of these systems is global. In this paper, we consider $M$ as a closed manifold. A $\mathsf{PSVF}$ defined on $M$ is a tangent vector field to $M$, which it’s just not differentiable in the points that belong to a submanifold $N$. The submanifold $N$ has the dimensions $n-1$ and, usually, is called discontinuity submanifold. Two approaches exist in the literature for formulating the equations for $\mathsf{PSVF}$; these are the equivalent control method of Utkin and the convex method of Filippov. In the book by  Filippov, see \cite{Filippov}, the author established conversions used in this paper to define solution orbit for a $\mathsf{PSVF}$. The problem of guaranteeing the existence of periodic orbits in smooth vector fields by parts is of colossal importance. In the literature, one of the tools studied for this purpose is the first recurrence map or Poincaré map; a fixed point this map corresponds to a periodic orbit. Some of the works in this context are: \cite{carmona}, \cite{du}, and \cite{llibre3}. In the papers \cite{llibre}, and \cite{llibre1}, the authors use the first integral to study the existence of crossing periodic orbits. Finally, the well-known theory of averaging was also used, for example, in the work \cite{llibre2}. Moreover, in  \cite{euzebio2014estudo}, \cite{junior2016orbitas} and \cite{tonon2010sistemas}, the author also study the existence of periodic orbits in piecewise smooth systems.
\newline

The Conley theory has been developed for continuous and discrete flows, multiflows, and semiflows. Referent to flows, one of the first works related to Conley Theory and discontinuous systems correspond to \cite{Casagrande}. The authors use a regularization of the discontinuous vector field, see \cite{SotomayorTeixeira1}, to adapt the Conley index for continuous flows. They define the D-Conley index and show it is invariant by homotopy. More recently, Cameron Thieme has been extending the Conley Theory for multiflows. His main objective is to generalize Conley index theory to differential inclusions having the Filippov systems as motivation. Firstly, in \cite{Cameron1}, he introduces differential inclusions and Filippov systems, and he shows the existence of a multiflow for this systems class. Subsequently in the preprints \cite{Cameron2} and \cite{Cameron3}, the author exposes a definition of perturbation (see  Definition 3.1 in \cite{Cameron2}) and shows that both the isolating neighborhoods and the attractor-repeller decomposition are stable, which are essential objects in the Conley theory. 
\newline

Our approach is to construct a semiflow for Filippov systems in order to apply the well-established results of the Conley theory. Recently, a related work has been done by Mrozek and Wanner in \cite{mrozek2020creating}; their engaging paper shows the construction of a continuous semiflow on a finite topological space $\mathsf{X}$ for a combinatorial vector field (a discretization of a piecewise smooth vector field).
\newline


This paper is structured as follows. The preliminary refer to Conley index and notation of $\mathsf{PSVF}$ is in Section II. Section III demonstrates the construction of a semiflow using the positive trajectories of a $\mathsf{PSVF}$ that, according to Filippov convention, only present crossing and sliding; there exists a unique solution for positive motion by the trajectories of the $\mathsf{PSVF}$. At the head of section III is the main result of this article. Finally, in Section IV, we are using the Conley Theory applicate to the semi-dynamical system construed in Section III to find periodic trajectories in $\mathsf{PSVF}$. 

\section{Preliminaries} Throughout this paper, we identify the intervals $(-\infty,0]$ and $[0,\infty)$ by $\mathbb{R^-}$ and $\mathbb{R^+}$, respectively. Let $\mathsf{X}$ be a topological space and $f$ a real-valued function on $\mathsf{X}$. We say that $f$ is upper semi-continuous at $x^*$ if and only if $\lim\sup_{{x}\to x^*} {f(x)}\leq f(x^*).$

\begin{definition}\label{semiflowdef}
The pair $(\mathsf{X},\phi)$ is called a (continuous) semi-dynamical system and $\phi$ a semiflow if $\mathsf{X}$ is a Hausdorff topological space and $\pi$ is a mapping, $\pi:\mathsf{X}\times \mathbb{R}^+\longrightarrow \mathsf{X}$ which satisfies
\begin{itemize}
\item [(1)] $\phi(x,0)=x$ for each $x\in \mathsf{X}$ (initial value property),
\item [(2)] $\phi(\phi(x,t),s)=\phi(x,t+s)$ for each  $x\in \mathsf{X}$ and
$t,s\in{\mathbb{R}^+}$ (semigroup property), and
\item [(3)] $\phi$ is continuous on the product space $\mathsf{X}\times\mathbb{R}^+$ (continuity property).
\end{itemize}
\end{definition}

Without loss of generality, if we are working with a fixed semiflow, we must often use the most straightforward notation  $x\cdot t$ in place of $\phi(x,t)$. A function $\sigma: I \longrightarrow \mathsf{X}$ where $I$ is a nonempty interval in $\mathbb{R}$ is called a solution of $(\mathsf{X},\phi)$ if whenever $t\in I$, $s\in{\mathbb{R}^+} $ and $t+s \in I$, then $\sigma(t)\cdot s=\phi(\sigma(t),s) = \sigma(t+s)$. The interval $I$ is the domain of $\sigma$ and according to our notation is represented by $D(\sigma)$. If $x\in \mathsf{X}$, a solution $\sigma$ with $0\in{D(\sigma)}$ and $\sigma(0)=x$ is called a solution through $x$. The function $\phi_x:\mathbb{R}^+\longrightarrow \mathsf{X}$ given by $\phi_x(t)=x\cdot t$ is a solution through $x$ and, indeed, is the unique solution through $x$ with domain $\mathbb{R}^+$. We call this solution the positive motion through $x$.

\subsection{Conley index}

In this section, we review the Conley index theory, following the definitions given by K. Mischaikow and M.  Mrozek \cite{Mischaikow}. For this section assume that $\mathsf{X}$ is a metric space and $\phi: \mathsf{X} \times \mathbb{R^+}\longrightarrow \mathsf{X}$ a semiflow. A subset $S \subset \mathsf{X}$ is said to be an \textit{invariant set} with respect to the semiflow $\phi$ if, for all $p\in S$, one has $p\cdot t \in S$ for all $t \in\mathbb{R^+}$. In other words, $\phi(S,\mathbb{R^+})=S$. Let $N\subset \mathsf{X}$ be a subset of $\mathsf{X}$. The \textit{maximal invariant set of $N$} is defined by: $\mathrm{inv}(N)=\{x\in{\mathsf{X}}| \ x\cdot t \in N, \text{ for all } t\in{\mathbb{R^+}}\}$. A subset $S\subset \mathsf{X}$ is called an \textit{isolated invariant set} if there exists a compact neighborhood $N$ of $S$ in $\mathsf{X}$ such that $S\subset \mathrm{int}(N)$ and $S=\mathrm{inv}(N)$. In this case, $N$ is said to be an \textit{isolating neighborhood} for $S$ in $\mathsf{X}$. The set $N$ is called an \textit{isolating neighborhood} for $\phi$ if it is closed, contained in the domain of $\phi$, and $\mathrm{Inv}(N)\subset \mathrm{int}(N)$. The most relevant property of the Conley index is its invariance by continuation, even under small perturbations.

\begin{definition}
	Let $S\subset \mathsf{X}$ be an isolated invariant set. A pair $(N,L)$ of compact sets in $\mathsf{X}$ is said to be an \textit{index pair} for $S$ in $\mathsf{X}$ if $L\subset N$ and
	\begin{enumerate}
		\item $\overline{N\setminus L}$ is an isolating neighborhood for $S$ in $\mathsf{X}$;
		\item $L$ is an positively invariant in $N$, that is, if $x\in{L}$ and $x\cdot [0,T]\subset N$ then $x\cdot[0,T]\subset L$;
		\item $L$ is the exit set of the semiflow, that is, if $x\in N$ and $x\cdot\mathbb{R^+}\varsubsetneq N $ then there exists $T>0$ such that $x\cdot[0,T]\subset N$ and $x\cdot T\in L$.
	\end{enumerate}
\end{definition}

\begin{figure}[H]
	\includegraphics[scale=0.8]{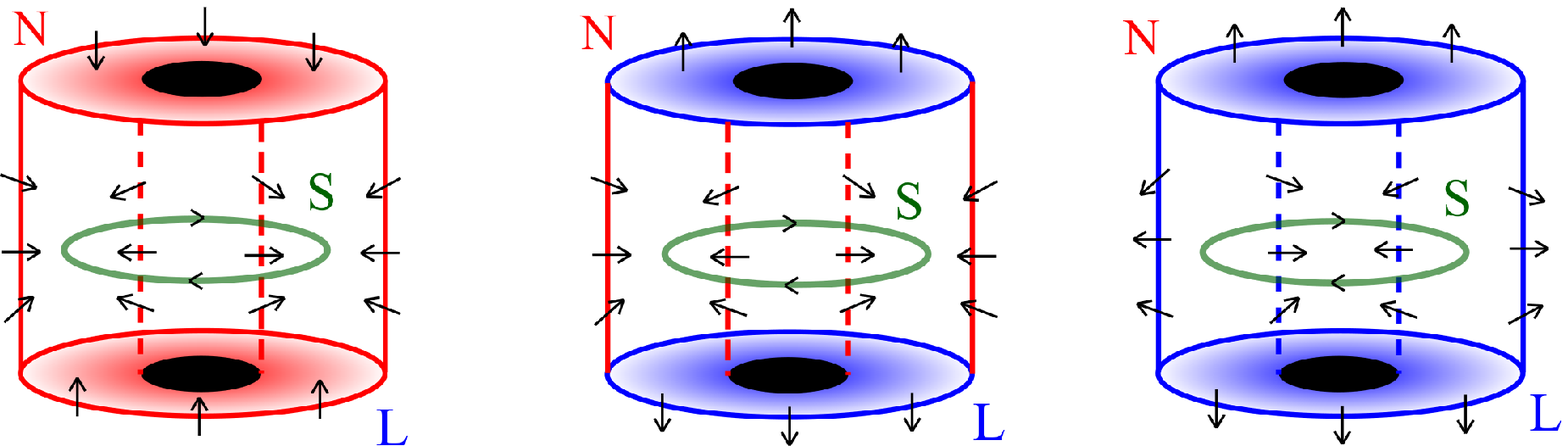}
	\caption{}
	\label{figure1}
\end{figure}

Given a pair $(N,L)$ of topological spaces with $L\subset N$ and $L\neq \emptyset $, define: 
\begin{equation}\label{relation}
x\sim y \Leftrightarrow x=y \text{ or  } x,y\in L.
\end{equation}
Denote by $N/L$ the pointed space $(N/\sim, [L])$. The Figure \ref{figure1} shows the index pair $(N,L)$ for a hyperbolic invariant set that is diffeomorphic to a circle in three dimensions.
Next, we present the definitions of the Conley index that we use throughout this paper.
\begin{definition}
	The \textit{Homotopy Conley Index} of $S$ is defined as the homotopy type of the pointed space $N/L$, where $(N,L)$ is the index pair $S$.
\end{definition}
Note that the definition that the homotopy Conley index is the homotopy type of a topological space. Unfortunately, operating with homotopy classes of spaces is extremely difficult. To evade this problem, it is useful to consider the cohomology Conley index.
\begin{definition}
Let $S$ be a isolated invariant set with respect to the semiflow $\phi$ and let $(N,L)$ be a index pair for $S$. The \textit{ cohomology Conley index} is defined as 
	\begin{equation}
	CH^*(S)=CH^*(S,\phi):=H^*(N/L)\approx H^*(N,L)
	\end{equation}
where $H^*$ denotes the Alexander-Spanier cohomology with integer coefficients.
\end{definition}

The most important property of an isolating neighborhood is that it is robust concerning perturbation. Still, before enunciating this property, the following definition is necessary.

\begin{definition}
Let $N\subset \mathsf{X}$ be a compact set. Let $S_\lambda=\mathrm{inv}(N,\phi_\lambda)$. Two isolated invariant sets $S_{\lambda_0}$ and $S_{\lambda_1}$ are related by continuation or $S_{\lambda_0}$ continues to $S_{\lambda_1}$ if $N$ is an isolating neighborhood for all $S_{\lambda}, \lambda\in{[\lambda_0,\lambda_1]}$.
\end{definition}

Now, we state the continuation theorem for the Conley index.

\begin{thm}[Continuation Property]
Let $S_{\lambda_0}$ and $S_{\lambda _0}$ be isolated invariant stes that are related by continuation. Then,
$$CH^*(S_{\lambda_0})\approx CH^*(S_{\lambda_1}).$$
\end{thm}
The proof of the Continuation Property can be found in \cite{Salamon}. Let us now consider, as an example,  the Conley index of a stable periodic orbit.

\begin{exam}\label{exam2.7}
The homotopy type of the pointed space $N/L$ of a stable periodic orbit in two dimensions is equal to the  homotopy type of $S^1 \vee S^0$ (see Figure \ref{figure2}), and so
\begin{equation*}
CH^k(S) \approx \left \{ \begin{matrix} \mathbb{Z} & \mbox{ }k=0, 1,
\\ 0 &  \mbox{ otherwise}.
\end{matrix}\right. 
\end{equation*}
\begin{figure}[H]
	\includegraphics[scale=1.1]{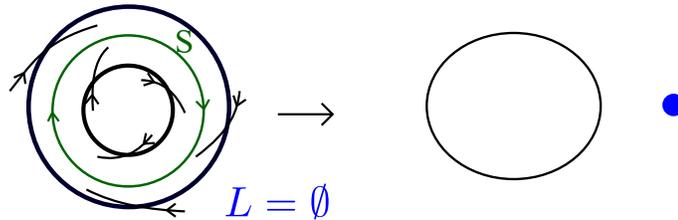}
	\caption{The homotopy type of the pointed space $N/L$ of a stable periodic orbit in dimension two.}
	\label{figure2}
\end{figure}	
\end{exam}

The following result generalizes Example \ref{exam2.7}.

\begin{figure}[H]
	\includegraphics[scale=0.9]{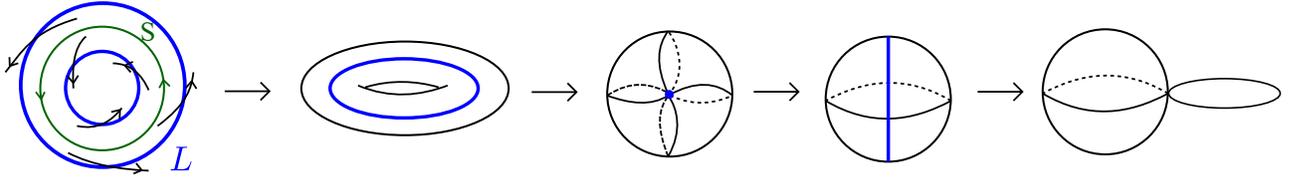}
	\caption{The homotopy type of a unstable periodic orbit in dimension two.}
	\label{figure3}
\end{figure}

\begin{prop}[Mischaikow, \cite{Mischaikow}]
Let $S$ be a hyperbolic periodic orbit with an oriented unstable
manifold of dimension $n+1$. Then
\begin{equation*}
CH^k(S) \approx \left \{ \begin{matrix} \mathbb{Z} & \mbox{ }k=n, n+1,
\\ 0 &  \mbox{ otherwise}.
\end{matrix}\right. 
\end{equation*}
\end{prop}
Figure \ref{figure3} shows the homotopy type of the pointed space $N/L$ of a unstable periodic orbit in two dimensions is equal to the homotopy type of $S^2 \vee S^1$.
Now, we refer to Poincaré's section that, the essential hypothesis to confirm the existence of periodic orbits in a semi-dynamic system when using the Conley theory. Suppose $N$ is an isolating neighborhood for a semiflow $\phi$. Then $\varXi \subset \mathsf{X}$ is a \textit{local section} for $\phi$ and $N$ if there exists $\xi> 0$ such that   $$C_{\xi}^N(\varXi):=\{x\in N | x\cdot(0,\xi)\cap\varXi\neq \emptyset \}$$ 
is open in $N$ and for every $x\in{C_{\xi}^N(cl(\varXi))}$ there exists a unique element of $x\cdot(0,\xi)\cap \varXi$.

\begin{definition}
	$\varXi \subset X$ is a \textit{Poincar\'{e} section} for $\phi$ in $N$ if:
	\begin{enumerate}
		\item $\varXi$ is a local section,
		\item $\varXi_N:=\{\varXi \cap N\}$ is closed and
		\item for every $x\in N$, $x\cdot(0,\infty)\cap \varXi \neq \emptyset.$
	\end{enumerate} 
\end{definition}

Observe that it is  not necessary to know $S$ to find a Poincaré section, also not required to be a subset of $N$. Certainly, if N has an exit set, then for any subset of N, there will be points in N whose orbits exit N before they cross again by the subset at issue, so no subset of N can be a Poincaré section. Suppose $(\Lambda,\rho)$ is a metric space that parameterizes a continuous family of semiflows $\phi^{\lambda}:\mathsf{X} \times \mathbb{R^+} \longrightarrow \mathsf{X}$. Furthermore, suppose that for some $\lambda_0\in{\Lambda}$, $N$ is an isolating neighborhood for $\phi^{\lambda_0}$ and $\varXi$ is a local section for $N$ and $\lambda_0$.

\begin{prop}\label{prop2.11}
Let $W$ be an admissible set for $\phi^{\lambda}$ for all $\lambda\in\Lambda$. Assume that $N\subset W$ is an admissible isolating neighborhood for $\phi^{\lambda}$ for all $\lambda\in\Delta\subset{\Lambda}$, a neighborhood of $\lambda_0$. Furthermore, assume that $N$ has a Poincaré section for $\phi^{\lambda_0}$ and that $\mathrm{inv}(N,\phi^{\lambda_0})=\emptyset$. Then, for $\Delta$ sufficiently small, given $\lambda\in{\Delta}$ there exists an isolating neighborhood $N_{\lambda}\subset W$ for $\phi^{\lambda}$ such that 
$$\mathrm{inv}(N_{\lambda},\phi^{\lambda})=\mathrm{inv}(N,\phi^{\lambda})$$
and $N_{\lambda}$ admits a Poincaré section.
\end{prop}

The previous proposition is the main result to guarantee the existence of periodic orbits that persists under perturbations of the system. In other words, for nearby semi-flows $\phi^{\lambda_0}$, there is an isolating neighborhood $N_{\lambda}$ that admits a section of Poincaré and has the same maximal invariant set for $\lambda\in{\Delta}$.
\newline

Given a Poincar\'{e} section $\varXi$ for $N$, there exists a subset $\varXi_0$ of $\varXi$, open in $\varXi$, such that $\varXi\cap S =\varXi_S\subset \varXi_0$ and such that, for every $x\in \varXi_0$, there exists a unique minimal strictly positive time $\pi_{\varXi}(x)$ with $x\cdot[0,\pi_{\varXi}(x)]\subset N$ and $x\cdot\pi_{\varXi}(x)\in\varXi$. The \textit{Poincar\'{e} map} $\varPi_{\varXi}$ associated with the Poincar\'{e} section $\varXi$ is defined by
\begin{equation}
\begin{aligned}
\varPi_{\varXi}: & \varXi_{0}\longrightarrow \varXi\\
& x \longmapsto x\cdot\pi_{\varXi}(x).
\end{aligned}
\end{equation}

A map $f:\mathsf{X} \longrightarrow \mathsf{Y}$ is \textit{locally compact} if every point $x\in{\mathsf{X}}$ admits a neighborhood $U$ such that the closure of $f(U)$ is compact. If, in addition, there exists a compact set $A$ such that for every $x\in{\mathsf{X}}$, $cl\{f^n(x)|n\in{\mathbb{N}}\}\cap A\neq \emptyset$, then $f$ is called a map of \textit{compact attraction}. A semiflow $\phi$ on $\mathsf{X}$ is  locally compact if there exists $t>0$ such that $\phi_t$ is locally compact and is of compact attraction if $\phi_t$ is a map of compact attraction for some $t>0$.

\begin{thm}[Mrozek, \cite{Mrozek}]
	Assume $f:\mathsf{X} \longrightarrow \mathsf{X}$ is a map of compact attraction and $K$ is an isolated invariant set for $f$. Then, the Conley index of $K$ under $f$ is of finite type.
\end{thm}
At the end of this section, we have the theorem that provides conditions for the Conley index of an isolated invariant set to contain a periodic trajectory. 

\begin{thm}[Mccord, Mrozek and Mischaikow,\cite{Mccord}]\label{teo1}
Assume $\mathsf{X}$ is an absolute neighborhood retract and $\phi:\mathsf{X}\times [0,\infty)\rightarrow \mathsf{X}$ is a semiflow with compact attraction. If $N$ is a isolating neighborhood for $\phi$ which admits a Poincar\'{e} section $\varXi$ and either 
\begin{equation}\label{eq1}
dim \ CH^{2n}(N,\phi)=dim \ CH^{2n+1}(N,\phi) \ \ \ \ \ for \ \ \ n\in{\mathbb{Z^+}}
\end{equation}
or 
\begin{equation}
dim \ CH^{2n}(N,\phi)=dim \ CH^{2n-1}(N,\phi) \ \ \ \ \ for \ \ \ n\in{\mathbb{Z^+}},
\end{equation}
	
where not all the above dimensions are zero, then $\phi$ has a periodic trajectory in $N$.
\end{thm}

\begin{cor}
	Under the hypotheses of Theorem \ref{teo1}, if $N$ has the Conley index of a hyperbolic periodic orbit, then $\mathrm{inv}(N)$ contains a periodic orbit.
\end{cor}

\subsection{Piecewise-smooth vector fields}
In this section, we introduce the Filippov convention for smooth piecewise vector fields, the reference \cite{Filippov}. Let $M$ be a closed $n$-dimensional $C^r$ manifold, denote by $\mathfrak{X}(M)$ the space of $C^r$ vector field tangent to $M$, endowed with the $C^r$-topology where $r$ is finite and sufficiently large. Let $\Sigma$ be a one codimension compact submanifold of $M$, that divides $M$ in two pieces, i.e.  $M=\Sigma^+\cup \Sigma^-$, where $\Sigma^+$ and $\Sigma^{-}$ are manifold with common boundary  $\partial \Sigma^+ =\partial \Sigma^-=\Sigma$. If $h: M \rightarrow \mathbb{R}$ is such that $h^{-1}(0)=\Sigma$, $h^{-1}([0,\infty))=\Sigma^+$, $h^{-1}((-\infty,0])=\Sigma^-$ and $0$  is a regular value $h$ so define a $\mathsf{PSVF}$ as follows
\begin{equation}\label{spvf}
Z(p) = \left\{ \begin{array}{rcl}
X(p)&& \text{if } p\in{\Sigma^{+}}, \\
Y(p) &&  \text{if } p\in{\Sigma^{-}}. 
\end{array}\right.
\end{equation}

We make the identification $Z=(X,Y)\in{\mathfrak{X}(M,h)}=\mathfrak{X}(M)\times \mathfrak{X}(M)$ equip with the product topology. We denote $\Sigma^\pm \setminus \Sigma$ for $\textrm{int}(\Sigma^\pm)$. 
\newline

For a vector field, $X\in{\mathfrak{X}(M)}$ the Lie derivative is defined as an operator $\mathcal{L}_{X}:\mathfrak{X}(M)\rightarrow \mathfrak{X}(M)$, such that, the Lie derivative of a vector field $Y$  in the direction of the vector field $X$ is the vector field $\mathcal{L}_{X}(Y)$ defined by 

$$\mathcal{L}_{X}(Y):=\frac{d}{dt}\bigg| _{t=0}(D\varphi(x,t))^{-1}Y(\varphi{(x,t)}),$$
where $(t,x)\mapsto\varphi{(x,t)}$ is the local flow of $X$. Generally $\mathcal{L}_{X}(Y)$ is defined by $[X,Y]$, where $[\cdot,\cdot]$ corresponds to the Lie bracket. The derivative for a differentiable function $h$ in $M$  is $\mathcal{L}_{X}(h)=(Xh)(p)=\sum_{i} m_i(p)\dfrac{\partial h}{\partial x_i}(p)$, where $ \bigg\{ \dfrac{\partial}{\partial{x_i}}\bigg\} $ is one basis associate to of a parametrization $\mathbf{x}:\mathcal{U}\subset{\mathbb{R}^n}\rightarrow M $ and each $m_i:\mathcal{U}\rightarrow \mathbb{R}$ is a function in  $\mathcal{U}$. Then, we have $\mathcal{L}_{X}(h)=\langle X(p),\nabla h(p)\rangle$, where $h$ indicates the expression of $h$ in the parametrization $\mathbf{x}$. Indeed, we can consider Lie derivatives
\begin{equation*}
 Xh(p)=\langle X(p),\nabla h(p)\rangle \text{ and } X^{k+1}h(p)=XX^kh(p)=\langle X(p),\nabla X^k h(p)\rangle,
\end{equation*}
where $\langle \cdot,\cdot \rangle$ denote the Euclidean inner product in $\mathbb{R}^{n}$ and $h$ indicates the expression of $h$ in a parametrization $\mathbf{x}$ of $M$. 
\newline 

Following the Filippov convention in \cite{Filippov}, we distinguish the following regions in the discontinuity manifold $\Sigma$:

\begin{itemize}
	\item[$\diamond$] the crossing  region: $\Sigma^{c}=\{p\in{\Sigma};Xh(p)Yh(p)>0\}$, $\Sigma^{c+}=\{p\in{\Sigma};Xh(p)>0 \text{ and } Yh(p)>0\}$ and  $\Sigma^{c-}=\{p\in{\Sigma}; Xh(p)<0 \text{ and } Yh(p)<0\}$,
	\item[$\diamond$] the escaping region:  $\Sigma^{e}=\{p\in{\Sigma};Xh(p)> 0 \text{ and }Yh(p)<0\}$,
	\item[$\diamond$] the sliding region: $\Sigma^{s}=\{p\in{\Sigma};Xh(p)<0 \text{ and }Yh(p)>0\}$.
\end{itemize}

\begin{definition}\label{slide}
	The sliding vector field associated to $Z$ is the vector field tangent to $\Sigma^{s}$ and defined at
	$$Z^{s}(p)= \frac{1}{Yh(p)-Xh(p)}(Yh(p)X(p)-Xh(p)Y(p)).$$  
	If $p\in{\Sigma^s}$ then $p\in{\Sigma^e}$ for $-Z$ and then we can define the escaping vector field on $\Sigma^e$ associated to $Z$ by $Z^e=-(-Z)^s$. 
\end{definition}

Let $X\in{\mathfrak{X}(\Sigma^+)}$. Since $\Sigma^+$ is a manifold with boundary then if $p\in{\partial \Sigma^+}$ we have that: $p$ is a \textit{regular point} of $X$ when $Xh(p)\neq 0$ or $p$ is a \textit{singular point} of $X$ when $Xh(p)=0$. We say that $p$ is a \textit{singularity tangency} of $X$ if the orbit passing for $p$ is tangent to boundary of $\Sigma^+ $ in p, i.e., $Xh(p)=0$ e $X(p)\neq0$, moreover, we say that $p$ is a \textit{singularity tangency} of $X$ if:

\begin{itemize}
\item[$\diamond$] $Xh(p)= 0$ and $X^2h(p)\neq 0$ (fold or singularity of order $2$),
\item[$\diamond$] $Xh(p)= X^2h(p)=0$, $X^3h(p)\neq 0$ and the set $\{Dh(p),DXh(p),  DX^2h(p)\}$ is linearly independent (cusp or singularity tangency of order $3$),
\item[$\diamond$] $\cdots$
\item[$\diamond$] $Xh(p)=X^2h(p)=\cdots= X^{m-1}h(p)=0$, $X^{m}h(p)\neq 0$ and the set $\{Dh(p),DXh(p), \\ DX^2h(p), \cdots,DX^{m-1}h(p)\}$ is linearly independent (singularity tangency of order $m$).
\end{itemize}

Although the Conley theory can be applied to semi-dynamical systems defined in any dimension, in this paper, we consider $M$ of dimension $3$, since it is the lowest dimension in which we find diversity in regions and tangencies formed in the discontinuity manifold.

\begin{definition}
	 Denoted by $S_X$ and $S_Y$ the tangency sets of $X$ and $Y$, respectively. If $Z=(X,Y)\in{\mathfrak{X}(M,h)}$, then the tangency set of $Z$ is given by $S_{Z}= S_{X} \cup S_{Y} $.  
\end{definition}

Denoted by $\Sigma_X=\{p\in{\Sigma}; \ Xh(p)\neq 0 \}$, respectively for $Y$ and $\Sigma_{Z}= \Sigma_{X} \cap \Sigma_{Y} $. 

\begin{rem} \label{rem2.19}
	Following \cite{Sotomayor}, $\Sigma$ is the disjoint union of submanifolds of decreasing dimension, if $S_X^j$ is the set of tangential singularities of order $j$ of $X$, then $\Sigma=\Sigma_X\cup S_X^2 \cup S_X^3$. This implies that every orbit meets the boundary of $\Sigma^+$ in a discrete set of points. Generically, a fold point of $X$ belongs to a local curve of fold points of $X$ with the same visibility, and cusp points occur as isolated points located at the extreme of curves of fold points. 
\end{rem}

The critical points of $\Sigma^s$ are consider \textit{pseudo-equilibria of $Z$} and the sliding vector field can be extended beyond the boundary of $\Sigma^{s}$, in fact, if $p\in{\partial \Sigma^s}$ then:
\begin{enumerate}
	\item if $Xh(p)=0$ but $Yh(p)\neq 0$ then  $Z^s(p)=X(p)$ ,
	\item if $Yh(p)=0$ but $Xh(p)\neq 0$ then  $Z^s(p)=Y(p)$  and
	\item if $Xh(p)=0$ and $Yh(p)=0$ then $p$ is a pseudo-equilibrium for $Z^s$ that is $Z^s(p)=p$.  
\end{enumerate}
Denote by  $\Sigma^S= \Sigma^s \cup \partial\Sigma^s$.
\newline 

If $p\in{\Sigma^c}$, then the orbit of  $Z=(X,Y)\in{\mathfrak{X}(M,h)}$ at $p$ is defined as the concatenation of the orbits of $X$ and $Y$ at $p$. Nevertheless, if $p\in{\Sigma\setminus\Sigma^c}$, there may be a lack of uniqueness of solutions. In this case, the flow of $Z$ is multivalued and any trajectory passing through $p$ originated by the orbits of $X$, $Y$, $Z^{s}$  and $Z^{e}$ is considered as a solution of $Z$. More details can be found in \cite{Filippov}. 
\newline

For $ \mathrm{dim}(M)=n$, on the discontinuity manifold, we assume only isolated crossing regions. For $ \mathrm{dim}(M) = 3 $, on the discontinuity manifold, it may have regions of crossing, sliding, and escaping, as well regions formed by regular, singular, and tangency points of order $ 2 $ and $ 3 $.  Given that we are interested in constructing the semiflow generated by trajectories of a $\mathsf{PSVF}$, in this paper, we do not consider systems with sliding and scaping regions simultaneously. For systems that have escaping and crossing regions, we consider the backward flow in order to construct a semi-flow.

\begin{definition}
	 Let $Z=(X,Y)\in{\mathfrak{X}(M,h)}$ a point $p\in{\Sigma}$ is said to be a \textit{ $\Sigma$-singularity} of $Z$ if $p$ is a tangential singularity, or a pseudo-equilibrium of $Z$. Otherwise, it is said to be a \textit{regular-regular} point of $Z$. 
\end{definition}

We said that a point $p\in{\Sigma}$ is said to be a \textit{fold point} of $X\in{\mathfrak{X}(\Sigma^+)}$ if $Xh(p)=0$ and  $X^2h(p)\neq 0$. If $X^2h(p)>0$ (resp. $X^2h(p)< 0$), then $p$ is a \textit{visible fold or with visible contact (resp. invisible fold or with invisible contact)}. Furthermore, a point $p\in{\Sigma}$ is said to be a \textit{cusp point} of $X\in{\mathfrak{X}(\Sigma^+)}$ if $Xh(p)=0$, $X^2h(p)=0$,  $X^3h(p)\neq 0$ and $\{Dh(p),DXh(p), DX^2h(p)\}$ is linearly independent. If $X^3h(p)>0$ (resp. $X^3h(p)< 0$), then $p$ is a \textit{visible cusp or with visible contact (resp. invisible cusp or with invisible contact)}. If  $Y\in{\mathfrak{X}(\Sigma^-)}$, the visibility condition is switched.
\newline

Resuming the previous results, we have the following definition for the local trajectory of one point $p\in M$ in a $\mathsf{PSVF}$.

\begin{definition}\label{traj}
The local trajectory $\varphi_{Z}(p,t)$ of a $\mathsf{PSVF}$ of the form (\ref{spvf}) through a point $p\in M$ is defined as follows where every interval $I$ contains zero.

\begin{enumerate}

\item  For $p\in{\textrm{int}(\Sigma^+)}$ and $p\in{\textrm{int}(\Sigma^-)}$ the trajectory is given by $\varphi_{Z}(p,t)= \varphi_{X}(p,t)$ and $\varphi_{Z}(p,t)= \varphi_{Y}(p,t)$ respectively, where $t\in{I\subset\mathbb{R}}$.

\item  For $p\in{\Sigma^{c+}}$ and taking the origin of time at $p$, the trajectory is defined as  $\varphi_{Z}(p,t)= \varphi_{Y}(p,t)$ for $t\in{I\cap \mathbb{R}^-}$ and  $\varphi_{Z}(p,t)= \varphi_{X}(p,t)$ for $t\in{I\cap \mathbb{R}^+}$. For the case $p\in{\Sigma^{c-}}$ the definition is the same reversing.

\item  For $p\in{\Sigma^{s}}$ and taking the origin of time at $p$, the trajectory is defined as $\varphi_{Z}(p,t)=\varphi_{Z^s}(t, p)$ for $t\in{I\cap \mathbb{R}^+}$ and $\varphi_{Z}(p,t)= \varphi_{1}(t, p)$ and $\varphi_{Z_1}(p,t)$ is either $\varphi_{X}(p,t)$, $\varphi_{Y}(p,t)$ or $\varphi_{Z^s}(p,t)$ for  $t\in{I\cap \mathbb{R}^-}$.

\item  For $p$ a generic singularity tangency and taking the origin of time at $p$, this is visible or invisible tangency for at least one of the fields $X$ ou $Y$, the trajectory is defined as $\varphi_{Z}(p,t)= \varphi_{1}(p,t)$ for  $t\in{I\cap\mathbb{R}^-}$ and $\varphi_{Z}(p,t)= \varphi_{2}(p,t)$ for $t\in{I\cap \mathbb{R}^+}$ where each $\varphi_{1}$, $\varphi_{2}$ is either $\varphi_{X}$ or $\varphi_{Y}$ or $\varphi_{Z^{\Sigma}}$. 

\item  For $p$ a equilibrium point, i.e., the equilibrium points of $X$, $Y$ and $Z^s$, the trajectory is defined as $\varphi_{Z}(p,t)=p$ for all $t\in \mathbb{R}$.

\end{enumerate}
\end{definition}

\section{Main Results}

The main objective of this section is to generate a semi-dynamical system with the forward trajectories of a piecewise-smooth vector field using the convex method of Filippov. Using this semiflow be used together with Theorem \ref{teo1}, we guarantee the existence of periodic orbits in piecewise-smooth vector fields without sliding region.
\newline

Let topological space $\mathsf{X}$ be a finite union of closed subsets $\mathsf{X}_i$, i.e.  $\mathsf{X}=\bigcup_{i=1}^n{\mathsf{X}_i}$; if for some topological space $\mathsf{Y}$, there are continuous maps $f_i:\mathsf{X}_i\longrightarrow \mathsf{Y}$ that agree on overlaps (i.e., ${f_i}{\mkern 1mu \vrule height 2ex\mkern2mu_{\mathsf{X}_i \cap \mathsf{X}_j}}={f_j}{\mkern 1mu \vrule height 2ex\mkern2mu_{\mathsf{X}_i \cap \mathsf{X}_j}}$), then the generalized pasting or gluing lemma says that there exists a unique continuous function $f:\mathsf{X} \longrightarrow \mathsf{Y}$ with $f\mkern 1mu \vrule height 2ex\mkern2mu_{\mathsf{X}_i}=f_i$, for all $i$. We used the gluing lemma with the functions of semiflows generated by positives trajectories of the vector fields $X$, $Y$, and $Z^s$. 
\newline

Throughout this paper, the dimension of the manifold $M$ is restricted to $n$-dimensional, $n \leq 3$. This assumption is needed to simplify the construction of a semiflow for the cases where the switching manifold $\Sigma$ has tangencies. For $n \geq 4$, the ideas will be similar; however, the challenge will be the analyze  of more tangency cases and defining Filippov convention for submanifolds with only tangencies. It is important to note that, for a Filippov system with only sliding and crossing regions, the assumption on the 3-dimensional manifold can be dropped, see Theorem \ref{teo3.23}.
\newline

Assume that $M$ is a closed $3$-dimensional $C^1$ manifold, and the discontinuity manifold has crossing regions and sliding regions. $M=\Sigma^+\cup \Sigma^-$, where $\Sigma^+$ and $\Sigma^{-}$ are manifold with common boundary $\Sigma$. We begin by introducing up some notations for subsets of $\Sigma$ that we use frequently.

\begin{definition}
For  $X\in{\mathfrak{X}(\Sigma^+)}$, denote $\Sigma_X^+=\{p\in\Sigma; \ Xh(p)>0\}$ and $\Sigma_X^-=\{p\in\Sigma; \ Xh(p)<0\}$ subsets of $\Sigma_X$. Furthermore,
	\begin{enumerate}
		\item $S_X^{v}=\{p\in{S_X}; \ p\text{ is a visible fold of } X\}$, 
		\item $S_X^{i}=\{p\in{S_X}; \ p\text{ is a invisible fold of } X\}$,
		\item $S_X^{ic}=\{p\in{S_X}; \ p\text{ is a cusp of } X \text{ and } X^3h(p)<0 \}$  and 
		\item $S_X^{vc}=\{p\in{S_X}; \ p\text{ is a cusp of } X \text{ and } X^3h(p)>0\}$.
	\end{enumerate}
Analogously to $Y\in{\mathfrak{X}(\Sigma^-)}$ and $Z^s\in{\mathfrak{X}(\Sigma^S)}$.
\end{definition}

In the following definition, $t_X^+(p)$ is the time necessary for the positive trajectory of $p$ to leave definitely $\Sigma^+$.  It is the primary tool when defining the domain in the semiflow. 

\begin{definition}\label{tx}
	Let $M$ be a closed manifold and $\Sigma$ ($\Sigma=h^{-1}(0)$ with $h: M \rightarrow \mathbb{R}$) a compact codimension $1$ submanifold of $M$, such that divides $M$ in two connected manifold $\Sigma^+$ and $\Sigma^{-}$ with common boundary $\Sigma$. For $X \in{\mathfrak{X}(M)}$, let $\Lambda_X^{+}=\{p\in{\Sigma^+}; \ \varphi_X(p,[0,\infty))  \nsubseteq  \Sigma^+\}$ and $t_X^+:\Sigma^+\longrightarrow \mathbb{R}^+\cup\{\infty\}$ such that
 		\begin{equation*}
		t_X^+(p) = \left\{ \begin{array}{lll}
			\ \infty     &&  \text{if }  p\notin{\Lambda_X^+},\\	
			\mathrm{inf} \{t>0; \ \varphi_X(p,t)\in{S_X^{ic} \cup \Sigma_X^-} \}
			&&  \text{if } \varphi_X(p,[0,t])\subset{\Sigma^+} \text{ for } t\in{I\cap\mathbb{R}^+} \text{ with } t\neq 0,\\
			\ 0  &&  \text{otherwise }  .
		\end{array}\right.
		\end{equation*}
	Analogously for $Y\in{\mathfrak{X}(M)}$.
\end{definition}

To better understand the belong definition, we present the following example with different cases to $t_X^+$.  

\begin{figure}[H]
	\includegraphics[scale=1.3]{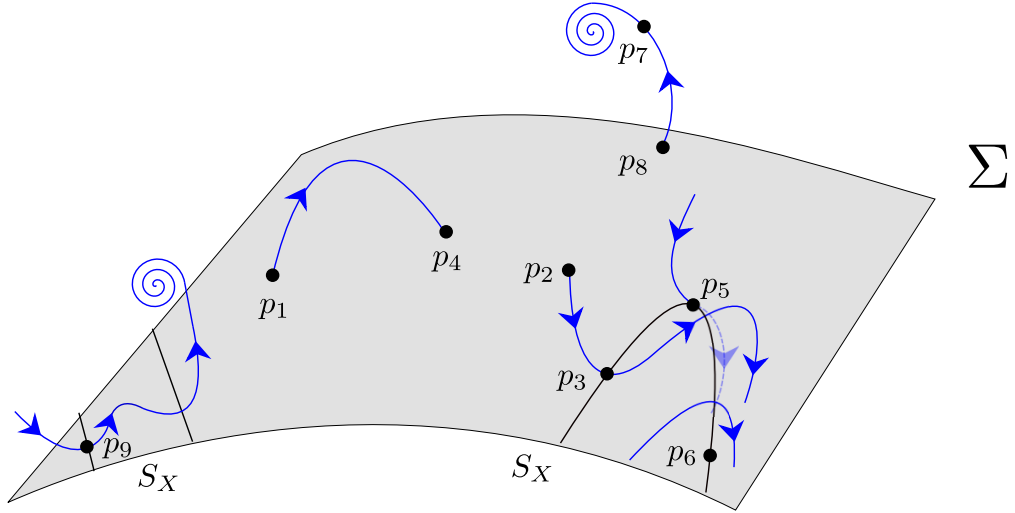}
	\caption{${p_i}\in{\Lambda_X^+}$ for $i=1,...,6$ but ${p_j}\notin{\Lambda_X^+}$ for $j=7,...,10$.}
	\label{figure4}
\end{figure}

\begin{exam}
	For $X\in{\mathfrak{X}(\Sigma^+)}$ we can have some points as in the Figure \ref{figure4}, note that ${p_i}\in{\Lambda_X^+}$ for $i=1,...,6$ but only for $p_1,p_2$ and $p_3$, $\varphi_X(p,[0,t])\subset{\Sigma^+}$ for $t\in{I\cap\mathbb{R}^+}$ with $t\neq 0$ then $t_X^+(p_i)>0$ for $i=1,2,3$ with $t_X^+(p_2)>t_X^+(p_3)$ and $t_X^+(p_i)=0$ for $i=4,5,6$. Now, for $p_i$ for $i=7,...,9$ we have that $\varphi_X(p,[0,\infty)) \subset \Sigma^+$ then, these points are not in $\Lambda_X^+$ and so $t_X^+(p_i)=\infty$ for $i=7,...,9$.
\end{exam}

\begin{figure}[H]
	\includegraphics[scale=1]{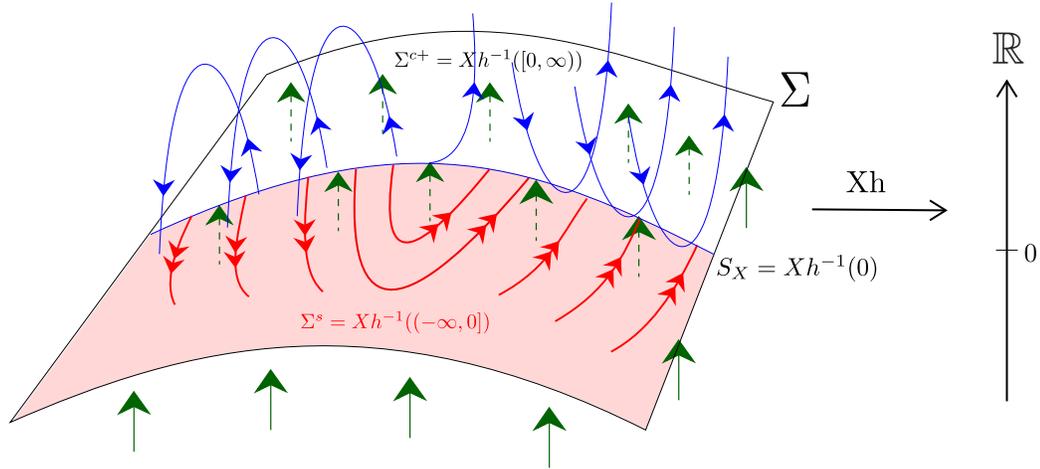}
	\caption{Local behavior of a cusp-regular singularities in a $\mathsf{PSVF}$.}
	\label{figure5}
\end{figure}

In the case, for the vector field $Z^s$, locally $\Sigma^s$ is a sub-manifold with boundary $\partial \Sigma^s=Xh^{-1}(0)$ or  $\partial \Sigma^s=Yh^{-1}(0)$. Without loss of generality, we take $p\in{S_X}$ such that $Yh(p)>0$ then locally the set $S_X$ that is a sub-manifold of co-dimension $1$ separates the discontinuity manifold into two pieces $\Sigma^S$ and $\Sigma^{c^+}$ such that $\Sigma^S=Xh^{-1}((-\infty, 0])$ and $\Sigma^{c+}=Xh^{-1}((0,\infty))$. See Figure \ref{figure5}.

In order to extend the Definition \ref{tx} to the vector field $Z^s$, in \cite{Teixeira1}, the following result is proved.

\begin{lem}\label{lem3.4}
\begin{itemize}
		\item [(i)]	The vector field $Z^s$ can be smoothly extended beyond the boundary of $\Sigma^s$. 
		\item [(ii)] If a point p in $\partial \Sigma^s$ is a fold point (resp. cusp point) of $X$ and a regular point of $Y$ then $Z^s$ is transverse to $\partial\Sigma^s$ at $p$ (resp. $Z^s$ has a quadratic contact with $\partial\Sigma^s$ at $p$).
\end{itemize}
\end{lem}
So, by Lemma \ref{lem3.4} and the following crollary, the Definition \ref{tx} is valid to $Z^s \in{\mathfrak{X}(\Sigma^S)}$.

\begin{cor}\label{cor3.5}
	Let $Z=(X,Y)$, $p\in{S_X}$ and $Yh(p)>0$. If $p\in{S_X^i}$ then $(Z^s)Xh(p)<0$ but if $p\in{S_X^v}$ then $(Z^s)Xh(p)>0$. If $p\in{S_X^{ic}}$ then $p\in{S_{Z^s}^v}$ and if $p\in{S_X^{vc}}$ then $p\in{S_{Z^s}^i}$ .
\end{cor}

In the following proposition, we show that the function $t_X^+$ is upper semi-continuous. Consider the notation, if $U$ is a open of $M$ then $\gamma_X(p)=\{\varphi_X(p,t); \ t\in{I}\}$ denoted the local orbit of a point $p\in{U}$;  $\gamma_X^+(p)$ and $\gamma_X^-(p)$ the local orbit positive and negative of $p$, respectively.

\begin{prop}\label{propo3.5}
	The application $t_X^{+}:\Sigma^+\rightarrow \mathbb{R}^+\cup{\{\infty\}}$ such that $p \mapsto t_X^{+}(p)$ is upper semi-continuous. Analogously for $t_{Z^s}^+$ and $t_Y^+$. 
\end{prop}

\begin{proof}
	Note that, if $p\notin{\Lambda_X^+}$ then $t_X^+(p)=\infty$ and $\lim\sup_{q\to p} t_X^+(q)\leq t_X^+(p)$ and so $t_X^{+}$ is upper semi-continuous at $p$. Now, assume that $p\in{\Lambda_X^+}$ with $p\in{\Sigma}$, let $q=\varphi_X(p,t_X^+(p))$ and $\gamma=\varphi_X(p,[0,t_X^+(p)])$. Choose coordinates $X=(x_1,x_2,x_3)$ around $p\in{\Sigma}$ such that $X=(1,0,0)$, let $x_3=g(x_1,x_2)$ be a $C^{\infty}$ solution of $h(x_1,x_2,x_3)=0$ with $g(0,0)=0$. Fix $N=\{x_1=0\}$ as the section transverse to $X$ at $p$. Define the projection $\sigma_{X_p}:(\Sigma,p)\longrightarrow (N,p)$ of $\Sigma$ along the orbits of $X$, onto $N$ is given by 
	$$\sigma_{X_p}(x_1,x_2,g(x_1,x_2))=(0,x_2,g(x_1,x_2)).$$
	
	Since $\Sigma$ is the disjoint union of submanifolds of decreasing dimension: $\Sigma_X\cup S_X^2 \cup S_X^3$. The map $\sigma_{X_p}$ is an immersion at $p$. This implies that every orbit meets the boundary of $\Sigma^+$ in a discrete set of points. Let $\U$ open of $p$ in $M$ and $\U_{\Sigma^+}=\U\cap \Sigma^+$.
	\begin{enumerate}
		\item Assume that $p\in{\Sigma_X}$.
		 \begin{enumerate}
		 	\item[(1.1)] If  $p\in{\Sigma_X^-}$ then $t_X^+(p)=0$ and so  for all $\tilde{p}\in{\U_{\Sigma^+}}$ (Case (a) of Figure \ref{figure6}), we have that $t_X^+(\tilde{p})\longrightarrow t_X^+(p)$ whenever $\tilde{p}\longrightarrow p$.
		 	\item[(1.2)] If  $Xh(p)>0$ and assume that not exists internal tangencies in the arc trajectory $\gamma$ and since $\Sigma=\Sigma_X\cup S_X^2\cup S_X^3$ assume that $q\in{\Sigma_X^-}$, by long tubular flow theorem there exists a tubular flow $(F,f)$ of $X$ such that $F\supset \gamma$, taking $F$ so small that the vector field in the box induced by $f$ and $X$ is the constant field $f_*X=(1,0,0)$ and choose the same coordinates $X=(x_1,x_2,x_3)$ around $p\in{\Sigma}$ (Figure \ref{figure6}). Consequently, again, for all $\tilde{p}\in{\U_{\Sigma^+}}$ we have that $t_X^+(\tilde{p})\longrightarrow t_X^+(p)$ whenever $\tilde{p}\longrightarrow p$. 
		 	\begin{figure}[H]
		 		\includegraphics{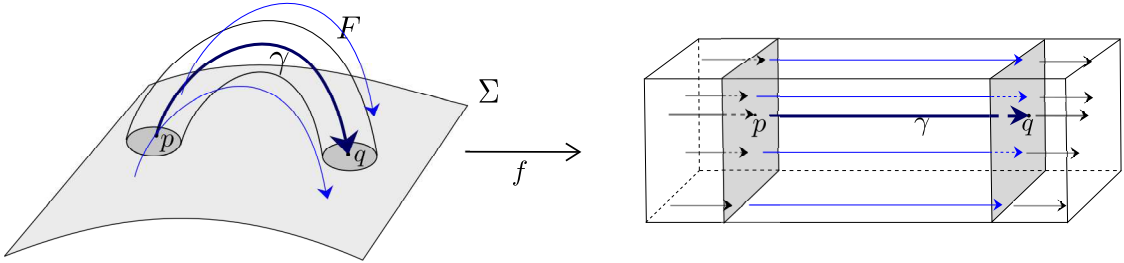}
		 		\caption{The tubular flow $(F,f)$ of $X$ such that $\gamma\subset F$.}
		 		\label{figure6}
		 	\end{figure}		 	
		 \end{enumerate}
		\item  Assume that $p$ is a singularity tangency of order two, then $p$ is a fold of $\sigma_X$ ($X^2h(p)=\frac{\partial^2 g }{\partial {x_ 1}^2}\neq 0$) and in this case there exits a $C^r$-diffeomorphism  $\tilde{\sigma}_{X_p}:(\Sigma,p)\longrightarrow (\Sigma,p)$, called the symmetric associated with $\sigma_{X_p}$ such that $\tilde{\sigma}_{X_p}(p)=p$, $\sigma_{x_p} \circ \tilde{\sigma}_{X_p} = \sigma_{X_p}$ and $\tilde{\sigma}_{X_p}^2=Id$. Observe that $S_X=\mathrm{Fix}(\tilde{\sigma}_X)$ and if $\tilde{p}\notin{S_X}$ then $\tilde{\sigma}_X(\tilde{p})$ is the point where the trajectory of $X$ passing through $\tilde{p}$ meets $\Sigma$.

		\begin{enumerate}
			\item[(2.1)] If $p\in{S_X^i}$ then $g(x_1,x_2)$ is conjugate to $(x_1,x_2)\longmapsto{x}_1^2$, $\tilde{p}=(\tilde{x_1},\tilde{x_2},\tilde{x_3})\in{\Sigma^+}$ iff $\tilde{x}_3\geq {\tilde{x}_1}^2$ (Case (b) of Figure \ref{figure7}) and so for all $\tilde{p}\longrightarrow p$ with $\tilde{p}\in{\Sigma^+}$ there exists a unique $t(\tilde{p})\geq 0$ such that the orbit-solution
			$t\longrightarrow \phi_X(\tilde{p},t)$ of $X$ through $\tilde{p}$ meets $\Sigma$, at a point $\tilde{q}=\phi_X(\tilde{p},t(\tilde{p}))$, so 
			
			$$\lim_{\tilde{p}\to p} t_X^+(\tilde{p})=\lim_{\tilde{p}\to p} t(\tilde{p})=t(p)=0=t_X^+(p).$$
			\item[(2.2)] Now, assume that $p\in{S_X^v}$ then $g(x_1,x_2)$ is conjugate to $(x_1,x_2)\longmapsto-{x_1}^2$, $\tilde{p}=(\tilde{x_1},\tilde{x_2},\tilde{x_3})\in{\Sigma^+}$ iff $\tilde{x}_3\geq -{\tilde{x}_1}^2$ (Case (c) of Figure \ref{figure7}). Reduce $\U$, if necessary, assume that $\gamma$ there are not internal tangencies in $\U$ then   $\U_{\Sigma^+}=\U_1\cup \U_2\cup \U_3$ with  $\U_1=\{\tilde{p}\in{\U_{\Sigma^+}}; \gamma_X^+(\tilde{p}) \cap \Sigma_X^- \neq \emptyset\}$, $\U_2=\{\tilde{p}\in{\U_{\Sigma^+}}; \gamma_X^+(\tilde{p})\cap S_X \neq \emptyset \text{ or } \gamma_X^+(\tilde{p})\cap \Sigma = \emptyset  \}$ and $\U_3=\{\tilde{p}\in{\U_{\Sigma^+}}; \gamma_X^-(\tilde{p})\cap \Sigma_X^+ \neq \emptyset\}$ then for the points $\tilde{p}\in{\U_1}$, $\gamma_X^+(\tilde{p}) \cap \Sigma_X^- \neq \emptyset$ and so $\lim_{\tilde{p}\to p} t_X^+(\tilde{p})=0$. Again, if we assume that $Xh(q)<0$ and there are no internal tangencies in the arc trajectory $\gamma$ so by long tubular flow theorem, reduce $F$, if necessary, for to not have internal tangencies in $F$ then for all $\tilde{p}\in{\U}$ there exists a unique $t(\tilde{p})\geq 0$ such that the orbit-solution $t\longrightarrow \phi_X(\tilde{p},t)$ of $X$ through $\tilde{p}$ meets $\Sigma$ (around of $q$) and for the points $\tilde{p} \in {\U_2 \cap \U_ 3}$ we have that $t(\tilde{p})=t_X^+(\tilde{p})$ so $\lim\sup_{\tilde{p}\to p} {t_X^+{(\tilde{p})}}=\inf\{\sup\{t_X^+(\tilde{p}); \ \tilde{p}\in{(\U_1\cup \U_2 \cup \U_3)\setminus{\{p\}}})\}\}\leq t_X^+(p).$
			\newline
			
			Now, assume that $Xh(q)<0$ but $\gamma$ there are, in fact finites, internal tangencies and let $p_1\in{\gamma}$ ($p_1\neq p$) the first point in $\gamma$ such that $p_1\in{S_X^v}$. Let $t_1>0$ such that $\varphi_X(p,t_1)=p_1$ then $t_1\ll t_X^+(p)$ and $\lim_{\tilde{p}\to p} t_X^+(\tilde{p})=t_1$ for the points $\tilde{p}\in{\U_3}$ and consequently $t_X^+$ is upper semi-continuous at $p$ in this cases.
					
		\end{enumerate}	
			\begin{figure}[H]
		\includegraphics{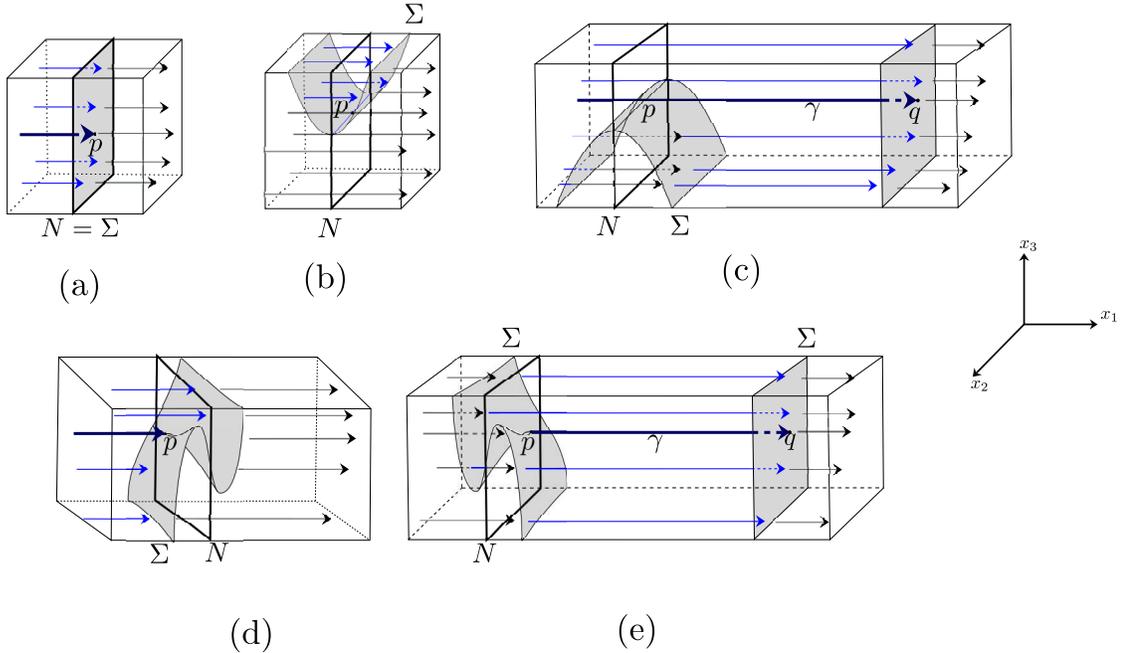}
		\caption{$X$ as a “straight” vector field and a “twisted” boundary.}
		\label{figure7}
		\end{figure}

		\item Assume that $p$ is a singularity tangency of order three, then $p$ is a cusp of $\sigma_X$.	
		
		\begin{enumerate}
		\item[(3.1)]
		If $X^3h(p)<0$ then $g(x_1,x_2)$ is conjugate to $(x_1,x_2)\longmapsto{x_1}^3+x_1x_2$, $\tilde{p}=(\tilde{x_1},\tilde{x_2},\tilde{x_3})\in{\Sigma^+}$ iff $\tilde{x}_3\geq {\tilde{x}}_1^3+\tilde{x}_1\tilde{x}_2$ (Case (d) of Figure \ref{figure7}) then $\gamma^+(\tilde{p})\cap \mathrm{int}({U \cap\Sigma^-}) \neq \empty$ and so $\lim_{\tilde{p}\to p} t_X^+(\tilde{p})=t_X^+(p)$. 
		
		\item[(3.2)]If $X^3h(p)>0$ then $g(x_1,x_2)$ is conjugate to $(x_1,x_2)\longmapsto-{x_1}^3-x_1x_2$ (Case (e) of Figure \ref{figure7})
		then $\U_{\Sigma^+}=\V_1\cup \V_2\cup \V_3$ with  $\V_1=\{\tilde{p}\in{\U_{\Sigma^+}}; \  \tilde{p}\in{S_X^i} \text{ or } \gamma_X^+(\tilde{p}) \cap \Sigma_X^- \neq \emptyset \}$, $\V_2=\{\tilde{p}\in{\U_{\Sigma^+}}; \ \gamma_X^+(\tilde{p})\in{S_X^v} \}$ and $\V_3=\{\tilde{p}\in{\U_{\Sigma^+}}; \ \tilde{p} \notin{ \V_1 \cup V_2} \}$. For the points in $\V_1$ by definition of $t_X^+$ we have that $\lim_{\tilde{p}\to p} t_X^+(\tilde{p})=0$, for the other points in $U_{\Sigma}$, by long tubular flow theorem there exists a tubular flow $(F,f)$ of $X$ such that $F\supset \gamma$ and consequently, to use Case $2.2$ above for the points in $\V_2 \cup \V_3$ for show that $\lim_{\tilde{p}\to p} t_X^+(\tilde{p})=t_X^+(\tilde{p})$, so $t_X^+$ is upper semi-continuous at $p$.
		
		\end{enumerate}			
	\end{enumerate}
 Now, if $p\in{\Lambda_X^+}$ but $p\notin{\Sigma}$ then we using the long tubular flow theorem for $\gamma$ (and using some the Cases 1), 2) or 3), if it is necessary) for show the assertion. And so $t_X^+$ is upper semi-continuous at for all $p\in{\Sigma^+}$ therefore it is a function upper semi-continuous.
\end{proof}

The following remark is a result notable of the before proposition. For each point in the manifold, there are only three possibilities for the limit superior of the function $t_X^+$.

\begin{rem}\label{cortx2}
	If $p\in{\Sigma^+}$ and  $\tilde{p}\longrightarrow p$ then $\displaystyle\lim\sup_{\tilde{p}\to p} {t_X^+{(\tilde{p})}}$ is $0$, $t_X^+(p)$ or the time $t>0$ such that $\varphi_X(p,t)$ is the first point of internal tangency of $\gamma$ depending of the cases $1)$, $2)$ or $3)$ of the proof of Proposition \ref{propo3.5}.
\end{rem}

A natural question is to ask whether function $t_X^+$ can be continuous. The answer is true in some cases. Using Proposition \ref{propo3.5} and Remark \ref{cortx2}  we have the following corollary.

\begin{cor}\label{cortx}
The application $t_X^{+}:\Sigma^+\rightarrow \mathbb{R}^+\cup{\{\infty\}}$ such that $p \mapsto t_X^{+}(p)$ is continuous whenever in $\Sigma=\Sigma_X\cup S_X^i$. Analogously for $t_{Z^s}^+$ and $t_Y^+$. 
\end{cor}

\begin{proof}
	Follow by long tubular flow theorem for the trajectory arc  $\gamma=\varphi_X(p,[0,t_X^+(p)])$ for $p\in{\Lambda_X^+}$.
\end{proof}

Now, we build the hypotheses to use the gluing lemma to create the semi-dynamical system. Again, we assume that $M$ is a closed $3$-dimensional $C^1$ manifold, and the discontinuity manifold has crossing regions and sliding regions. The boundary of these regions is composed of tangency points. We consider the points in $\Sigma$ classified in two types $\mathrm{A}$ and $\mathrm{B}$, as follows:

\begin{itemize}
	\item [Case A.] Singularity points of type singular-regular;
	
	\begin{itemize}
		\item [A1.]	$p\in{S_X^v}$ and $Yh(p)> 0$ (or $p\in{S_Y^v}$ and $Xh(p)< 0$);
		\item [A2.] $p\in{S_X^i}$ and $Yh(p)> 0$ (or $p\in{S_Y^i}$ and $Xh(p)<0$);
		\item [A3.]	$p\in{S_X^{ic}}$ and $Yh(p)> 0$ (or $p\in{S_Y^{ic}}$ and $Xh(p)<0$);
		\item [A4.]	$p\in{S_X^{vc}}$ and $Yh(p)> 0$ (or $p\in{S_Y^{vc}}$ and $Xh(p)<0$).
	\end{itemize}
	
	\begin{figure}[h]
		\includegraphics[scale=1.2]{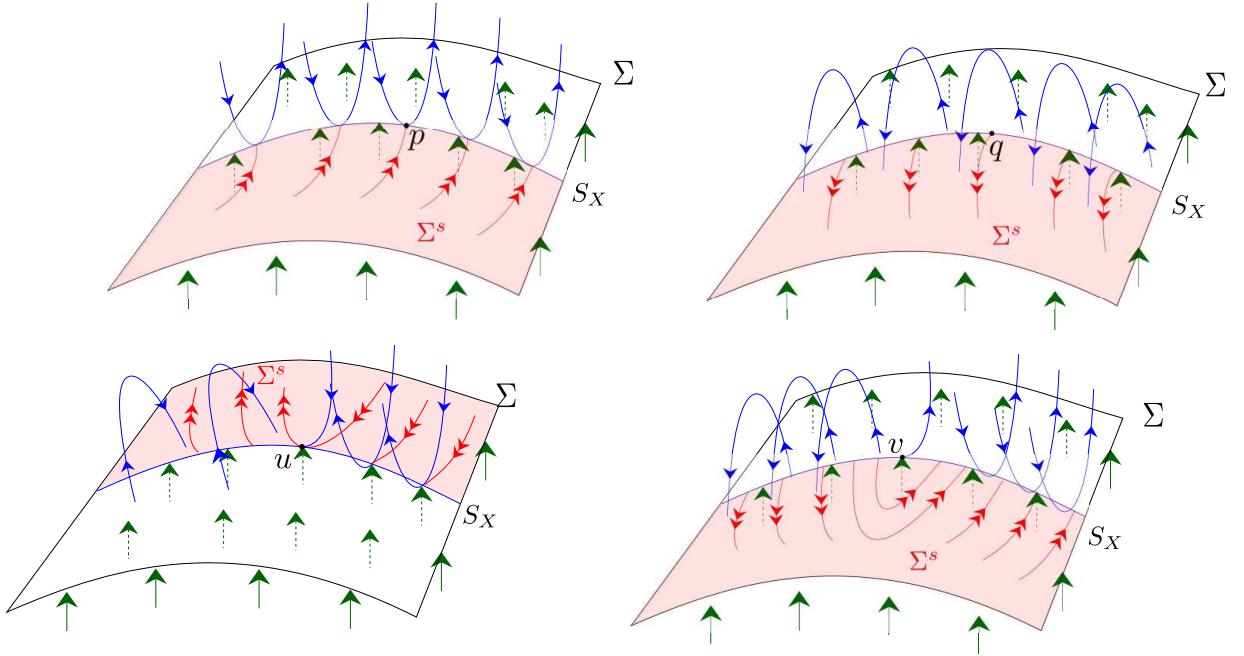}
		\caption{Local behavior of the points in a case $\mathrm{A}$.}
	\end{figure}
	\item [Case B.] Singularity points of type singular-singular;
	
	\begin{itemize}
		\item [B1.] $p\in{S_X^v \cap S_Y^i}$ (or $p\in{S_X^i \cap S_Y^v}$) with $p\in{\partial \Sigma^{c+}\cap \partial \Sigma^{c-}}$;
		\item [B2.] $p\in{S_X^{ic}}$ and $p\in{S_Y^i}$ (or $p\in{S_Y^{ic}}$ and $p\in{S_X^i}$) with $S_X$ and $S_Y$ tangent of $p$.
	\end{itemize}

\begin{figure}[h]
	\includegraphics[scale=1.2]{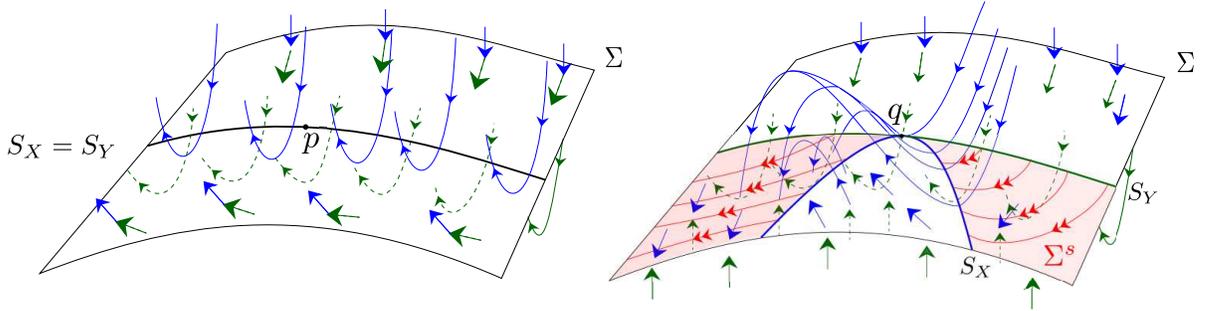}
	\caption{Local behavior of the points in a case $\mathrm{B}$.}
\end{figure}

\end{itemize}
We denoted, also by $A_1$ to $B_2$ the sets of points in the cases $A_1$ to $B_2$. Note that, by Definition \ref{slide}, if $p$ is the point of case  $\mathrm{B2}$ then $Xh(p)=Yh(p)=0$ and so $Z^sh(p)=0$. With the following lemma, we guarantee that the positive local trajectory for the tangency points in the cases $A1$ to $B2$ has uniqueness.

\begin{lem}\label{traj1}
	The positive local trajectory $\varphi_{Z}(p,t)$ of a $\mathsf{PSVF}$ of the form (\ref{spvf})  through a point $p\in S_Z$ in some of the cases $\mathrm{A1}$ to $\mathrm{B2}$ is unique.
\end{lem}
\begin{proof}
In fact, by Definition \ref{traj} and Lemma \ref{cor3.5}
	\begin{enumerate}
		\item [1)] for $p$ in  $\mathrm{A1}$,  $\mathrm{A4}$ and  $\mathrm{B1}$, the trajectory is given by $\varphi_{Z}(p,t)=\varphi_{X}(p,t)$ ( $\varphi_{Z}(p,t)=\varphi_{Y}(p,t)$);
		\item [2)] for $p$ in  $\mathrm{A2}$ and  $\mathrm{A3}$, the trajectory is given by $\varphi_{Z}(p,t)=\varphi_{Z^s}(p,t)$;
		\item [3)] for $p$ in  $\mathrm{B2}$, the trajectory is given by $\varphi_{Z}(p,t)=\varphi_{Z^s}(p,t)=p$
	\end{enumerate}
for $t\in{I\cap \mathbb{R^+}}$ and $0\in{I}$.
\end{proof}

In the following definition, for each $p\in{M} $ we build a sequence of points under the discontinuity manifold using the function $t_X^+$, also, together to one sequence of fields $X$, $Y$ and $Z^s$ and closed intervals of positive time. The objective of these sequences is to use them to restrict the domains of the flows generated by the fields in question.

\begin{definition}\label{omegax}
Let $Z=(X,Y)\in{\mathfrak{X}(M,h)}$ with $\Sigma=\Sigma^c \cup \Sigma^s\cup S_Z$ and $p\in S_Z$ in some of the cases $\mathrm{A1}$ to $\mathrm{B3}$. If $p\in{M}$ then denoted $p_0=p$ and $p_i=		\varphi_{Z_{{i-1}(p)}}(p_{i-1},t_{Z_{{i-1}(p)}}^+(p_{i-1}))$ for $i\in{\mathbb{Z}^+}$, where  $Z_{0(p)}=X$ if $p\in{\mathrm{int}(\Sigma^+)}$, $Z_{0(p)}=Y$ if $p\in{\mathrm{int}(\Sigma^-)}$ and
		
		
	\begin{equation*}
		Z_{i(p)} = \left\{ \begin{array}{rcl}
		X && \text{if }  p_{i}\in\Sigma^{c+} \text{or }  p_{i}\in{( \mathrm{A1} \cup \mathrm{A4} \cup \mathrm{B1})\cap(S_X^v \cup S_X^{vc}) },\\
		Z^s &&  \text{if } p_{i}\in\Sigma^{s} \text{or }  p_{i}\in {\mathrm{A2}\cup \mathrm{A3} \cup \mathrm{B2}},\\
		Y &&  \text{if } p_i\in{\Sigma^{c-}}  \text{or }  p_{i}\in{( \mathrm{A1} \cup \mathrm{A4} \cup \mathrm{B1})\cap(S_Y^v \cup S_Y^{vc}) }.
		\end{array}\right.
	\end{equation*}
		
Let $\Delta_X^+(p)=\{i\in{\mathbb{Z}^+}; \ Z_{i(p)}=X\}$ and $\Omega_X=\{(p,t)\in{M\times \mathbb{R}^+}; \ t\in{I_X(p)}\}$ whit $$I_X(p)=\bigcup_{i\in{\Delta_X^+}(p)} I_X^i(p)$$ 
where $I_X^0(p) = \left[ 0,a_{0(p)} \right] $,  $I_X^i(p) =	\left[ a_{{i-1}(p)},a_{i(p)} \right] $ for $i>0$ such that $a_{i(p)}=\displaystyle\sum_{j=0}^{i} t_{Z_{j(p)}}^+(p_j)$. Analogously, we defined $\Omega_Y$ and $\Omega_{Z^s}$ for the vector fields $Y$ and $Z^s$, respectively.
\end{definition}

Now, we define the functions to be glued together.

\begin{definition}\label{semiflow} Let $\varphi_X$, $\varphi_Y$ and $\varphi_{Z^s}$ the flows of the vector fields  $X\in{M}$, $Y\in{M}$ and $Z^s\in{\mathfrak{X}(\Sigma^S)}$, respectively. Denote by the restriction of the domains to $\Omega_X$,  $\Omega_Y$ or  $\Omega_{Z^s}$ correspond 
	$$\phi_{X}: \Omega_X \longrightarrow \Sigma^+,$$ 
	$$\phi_{Y}: \Omega_Y \longrightarrow \Sigma^- \text{ and }$$ 
	$$\phi_{Z^s}: \Omega_{Z^s} \longrightarrow \Sigma^S.$$ 
For $(p,t)\in{\Omega_X}$ then $t\in I_X^i(p)$ with $i\in{\Delta_X^+(p)}$ and $\phi_X(p,t)=\varphi_X(p_i,t-a_{i-1}(p))$. Analogously, for $\phi_Y$ and $\phi_{Z^s}$.
\end{definition}

\begin{rem}\label{rem3.14}
The $\phi_{X}: \Omega_X \longrightarrow M $ is a continuous function, in fact $\phi_X=\mathfrak{i}\circ \varphi_X|_{\Omega_X} $ with $\mathfrak{i}$ is the inclusion function of $\Sigma^+$ in $M$. Analogously, for $\phi_Y$ and $\phi_{Z^s}$.
\end{rem}
The following lemma is a significant result for the proof that the collage of the limits points of the domain of the functions for the gluing is correct. 
\begin{lem}\label{lem3.15}
If $(p,s)\in{\overline{\Omega_X}}$ then $\phi_{Z_{i(p)}}(p,s)\in{\Sigma^+}$ where $s\in{I_{Z_{i(p)}}^i(p)}$.
\end{lem}
\begin{proof}
Assume by contra-position that $q=\phi_{Z_{i(p)}}(p,s)\notin{\Sigma^+}$ then $q\in{\mathrm{int}(\Sigma^-)}$. Let $V$ a open neighborhood  of $\mathrm{int}(\Sigma^-)$ such that $q\in{V}$ then by Remark \ref{rem3.14} $\phi_Y^{-1}(V)$ is a open neighborhood of $M\times \mathbb{R^+}$ such that $(p,s)\in{\phi_Y^{-1}(V)}$ and for all $(\tilde{p},\tilde{s})\in{\phi_Y^{-1}(V)}$, $\tilde{s}\in{I_Y(p)\setminus I_X(p)}$ and so $\phi_Y^{-1}(V)\cap \Omega_X =\emptyset$ that is $ (p,s)\notin{\overline{\Omega_X}}$.
\end{proof}

For $(p,s)\in{M\times\mathbb{R^+}}$ let $\gamma_i(p)=\varphi_{Z_{i(p)}}(p_i,[0,t_{Z_{i(p)}}^+(p_i)])$ for each $i\in{\Delta_Z^+(p)}$ and $\Gamma$ the concatenation of for all $\gamma_i(p)$. Let $\Delta_Z^+(p)=\Delta_X^+(p)\cup\Delta_Y^+(p)\cup\Delta_{Z^s}^+(p)$.

\begin{prop}\label{prop3.15}
	Let $(p,s)\in\overline{\Omega_X}$ with $s\in{I_{Z_{i(p)}}^i(p)}$ and $i\in{\Delta_Z^+(p)}$ then $$\lim\limits_{\substack{\tilde{p} \to p \\ \tilde{s}\to s}}  \phi_X(\tilde{p},\tilde{s})=\phi_{Z_{i(p)}}(p,s)$$ whenever $(\tilde{p},\tilde{s})\in{(V\times I)\cap\Omega_X}$ with $V\times I $ is a neighborhood of $(p,s)$ in $M\times\mathbb{R^+}$.  Analogously for  $(p,s)$ in $\overline{\Omega_Y}$ and $\overline{\Omega_{Z^s}}$. 
\end{prop}

\begin{proof}
If $(p,s)\in{\Omega_X}$ follow by continuity of $\phi_X$. For $(p,s)\notin{\Omega_X}$ then $s\in{(I_Y(p)\cup I_{Z^s}(p))\setminus I_X(p)}$, that is, there exists $i\in{\Delta_Z(p)^+}$ such that $s\in{I_{Z_{i(p)}}^i(p)}$ with $Z_{i(p)} \neq  X $. Fixing $p\in{M}$, we proceed by induction on $i$.
\newline

\textit{Induction basis:} Assume that $i=0$ then $(p,s)\in{M\times\mathbb{R^+}}$ is such that $s\in{I_{Z_{0(p)}}^0(p)}$, note that, $p\in{\Sigma^-}$. Firstly, assume that $p\in{\mathrm{int}(\Sigma^-)}$ then $Z_{0(p)}=Y$ and by Lemma \ref{lem3.15} $s\neq 0$. Then $\gamma_0(p)=\varphi_Y(p,[0,a_{0(p)}])=\varphi_Y(p,[0,t_Y^+(p)])$ and so we have the following cases:
\begin{enumerate}
	\item For $s<a_{0(p)}$, this is $\gamma_0(p)$ has internal tangencies. The first internal tangency is a adherent point of $\Omega_X$ if, and only if, these is a point of case $\mathrm{B1}$ then let $q=\phi_Y(p,s)=\varphi_Y(p,s)$ this point, so $q\in{S_X^i\cap S_Y^v}$ and also $q\in{\partial \Sigma^{c+} \cap \partial \Sigma^{c-}}$. We are going to show that $\phi_X(\tilde{p},\tilde{s})\longrightarrow q$ when $\tilde{p}\longrightarrow p$ and $\tilde{s}\longrightarrow s$. Note that (Case (a) of Figure \ref{figure10}) $\tilde{s}\in{I_X^1(\tilde{p})}$ and as $p\in{\mathrm{int}(\Sigma^-)}$ and $\gamma_0(p)$ has internal tangencies by Remark \ref{cortx2} we have that $\displaystyle\lim_{\tilde{p}\to p}t_Y^+(\tilde{p})=s$ and $$\lim_{\tilde{p} \to p } \tilde{p}_1= \lim_{\tilde{p} \to p} \varphi_Y(\tilde{p},t_Y^+(\tilde{p}))=\varphi_Y(\lim_{\tilde{p}\to p} \tilde{p},\lim_{\tilde{p}\to p} t_Y^+(\tilde{p}))=\varphi_Y(p,s)=q,$$ 
and so
	\begin{align*}
	\lim\limits_{\substack{\tilde{p} \to p \\ \tilde{s}\to s}}  \phi_X(\tilde{p},\tilde{s}) &=\lim\limits_{\substack{\tilde{p} \to p \\ \tilde{s}\to s}}  \varphi_X(\tilde{p}_1,\tilde{s}-a_{0(\tilde{p})})
	\\ &=\lim\limits_{\substack{\tilde{p} \to p \\ \tilde{s}\to s}}  \varphi_X(\tilde{p}_1,\tilde{s}-t_Y^+(\tilde{p}))
	\\ &= \varphi_X \left( \lim_{\tilde{p}\to p} \tilde{p}_1,\lim_{\tilde{s}\to s}\tilde{s}-\lim_{\tilde{p}\to p} t_Y^+(\tilde{p})\right)
	\\ &= \varphi_X(q, s-s)=q.
	\end{align*}
 
\begin{minipage}{\linewidth}
	\begin{minipage}{0.47\linewidth}
		\begin{figure}[H]
			\includegraphics[scale=1.2]{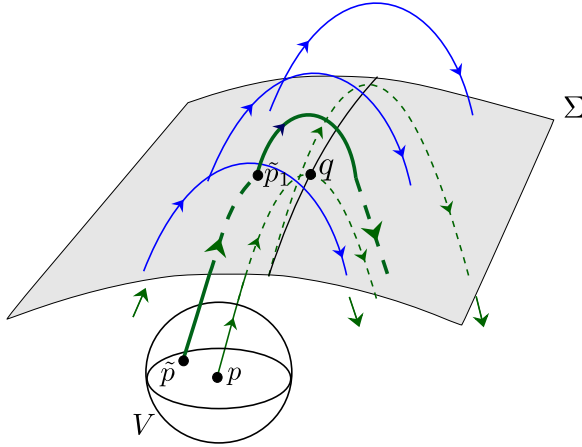}
			\caption{$q=\phi_Y(p,s)$ is a point of case $\mathrm{B1}$ and $\tilde{p}_1\in{\Sigma^{c+}}$.}
			\label{figure10}
		\end{figure}
	\end{minipage}
	\hspace{0.02\linewidth}
	\begin{minipage}{0.47\linewidth}
		\begin{figure}[H]
			\includegraphics[scale=1.2]{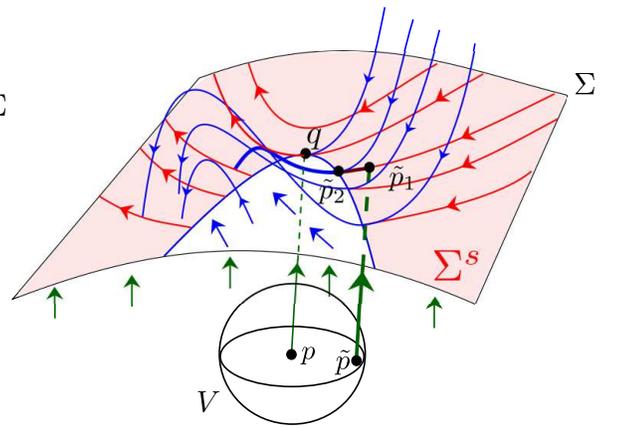}
			\caption{ $q=\phi_Y(p,s)$ is a point of case $\mathrm{A3}$, $\tilde{p}_1\in{\Sigma^{s}}$ and $\tilde{p}_2\in{\partial\Sigma^{s}}$.}
			\label{figure11}
		\end{figure}
	\end{minipage}
\end{minipage}		
\newline	

	Now, if $\varphi_Y(p,s)$ is not the first internal tangency, the proof follow by induction in $k$ where $k$ is the position of internal tangency, in fact $\gamma_0(p)$ has finite tangencies with $\Sigma$. 
	\newline
	
	\item For $s=a_{0(p)}$, assume that $\gamma_{0(p)}$ has no internal tangencies. In this case, $q=\phi_Y(p,s)=\phi_Y(p,a_{0(p)})=\phi_Y(p,t_Y^+(p)) =p_1$ is a point of case either $\mathrm{A1}$ or $\mathrm{A3}$. For $q \in{\mathrm{A1}}$ the proof is analogous to case (1), it remains to prove to $q\in{\mathrm{A3}}$. As $(p,s)$ is a adherent point of $\Omega_X$ then either $\tilde{s}\in{I_X^1(p)}$ if $\tilde{p}_1\in{\Sigma^{c+}\cup\partial\Sigma^s}$ (is analogous to case $(1)$) or $\tilde{s}\in{I_X^2(p)}$ if $\tilde{p}_1\in{\Sigma^s}$ (Case (b) of Figure \ref{figure11}). For $\tilde{s}\in{I_X^2(p)}$ then $\tilde{p}_1\longrightarrow q$ when $\tilde{p}\longrightarrow p$ and by using the part $2.2$ of the proof of Proposition \ref{propo3.5} but with the map $t_{Z^s}^+:\Sigma^S\longrightarrow \mathbb{R^+}$ we have that  $\displaystyle\lim_{\tilde{p}_1\to q}t_{Z^s}^+(\tilde{p}_1)=0$, in fact $\gamma_1(\tilde{p})=\varphi_{Z^s}(\tilde{p}_1,[0,t_{Z^s}^+(\tilde{p}_1)])$ has no internal tangencies with $\partial\Sigma^s$ and so
	
	$$\lim_{\tilde{p}\to p} \tilde{p}_2= \lim\limits_{\substack{\tilde{p} \to p \\ \tilde{p}_1\to q}} \varphi_{Z^s}(\tilde{p}_1,t_{Z^s}^+(\tilde{p}_1))=\varphi_{Z^s}(\lim_{\tilde{p}\to p} \tilde{p}_1,\lim_{\tilde{p}_1\to q} t_{Z^s}^+(\tilde{p}_1))=\varphi_{Z^s}(q,0)=q,$$ 
	and so
	\begin{align*}
	\lim\limits_{\substack{\tilde{p} \to p \\ \tilde{s}\to s}}  \phi_X(\tilde{p},\tilde{s}) &=\lim\limits_{\substack{\tilde{p} \to p \\ \tilde{s}\to s}}  \varphi_X(\tilde{p}_2,\tilde{s}-a_1(\tilde{p}))
	\\ &=\lim\limits_{\substack{\tilde{p} \to p \\ \tilde{s}\to s}}  \varphi_X(\tilde{p}_2,\tilde{s}-(t_Y^+(\tilde{p})+t_{Z^s}^+(\tilde{p}_1)))
	\\ &= \varphi_X(\lim_{\tilde{p}\to p} \tilde{p}_1,\lim_{\tilde{s}\to s}\tilde{s}-(\lim_{\tilde{p}\to p} t_Y^+(\tilde{p}) + \lim_{\tilde{p}_1\to q} t_{Z^s}^+(\tilde{p}_1)  ))
	\\ &= \varphi_X(q, s-(s+0))=q.
	\end{align*}	
	If $\gamma_{0(p)}$ has internal tangencies and $(p,a_{0(p)})\in\overline{\Omega_X}$ then this internal tangencies are  points of case $\mathrm{B1}$ and so use to item (1) and then the item (2) but instead of $\tilde{p}$ use to $\tilde{p}_j$ with $j\in{\Delta_Z^+{\tilde{p}}}$ such that $\tilde{p}_j \to \tilde{q}$ where $p^*$ is the last of the internal tangencies of $\gamma_{0(p)}$ such that  $(p,s^*)\in\overline{\Omega_X}$ with $\phi_{Z_{0(p)}}(p,s^*)=\phi_Y(p,s^*)=q^*$.
\end{enumerate}
Now, assume that $p\in{\Sigma}$. As $s\in{I_{Z_{0(p)}}^0(p)}$ then  
either $s=0$ or $0<s\leq t_{Z_{0(p)}}^+(p)$. If $s=0$ then $p$ is a point of cases $A2$, $A3$ or $B1$ and doing a similar analysis to the previous items $(1)$ and $(2)$ we have that $\lim\limits_{\substack{\tilde{p} \to p \\ \tilde{s}\to 0}}  \phi_X(\tilde{p},\tilde{s})=p$. If $0<s\leq t_{Z_{0(p)}}^+(p)$ then $\lim\limits_{\substack{\tilde{p} \to p \\ \tilde{s}\to s}}  \phi_X(\tilde{p},\tilde{s})$ is either a internal tangency of $\gamma_0(p)$ ( with $\Sigma$ when $Z_{0(p)}=Y$ or with $\partial\Sigma^s$ when $Z_{0(p)}=Z^s$ ) or $p_1$, and the proof is analogous to the previous items $(1)$ and $(2)$.
Analogously for  $(p,s)$ in $\overline{\Omega_Y}$ and $\overline{\Omega_{Z^s}}$ with $s\in{I_{Z_{0(p)}}^0(p)}$.
\newline 

\textit{Inductive step:} We assume that the proposition is true for all $(p,s)\in{\overline{\Omega_X}}$ ($(p,s)\in{\overline{\Omega_Y}}$ or $(p,s)\in{\overline{\Omega_{Z^s}}}$) with $s\in{I_{Z_{i(p)}}^i(p)=[a_{{i-1}(p)},a_{i(p)}]}$ thus for
$s=a_{i(p)}$, if $(p,a_{i(p)})\in{\overline{\Omega_{Z^*}}}$ with $Z^*=X,$ $Y$ or $Z^s$ so 
$$\lim\limits_{\substack{\tilde{p} \to p \\ \tilde{s}\to a_{i(p)}}}  \phi_{Z^*}(\tilde{p},\tilde{s})=\phi_{Z_{i(p)}}(p,a_{i(p)})=p_{i+1},$$
then there exists $j\in{\Delta_Z^+(\tilde{p})}$ such that $\tilde{p}_j\longrightarrow p_{i+1}$ and by continuity of $\varphi_{Z_{{j-1}(\tilde{p})}}$; $a_{{j-1}(\tilde{p})}\longrightarrow a_{i(p)}$ whenever $\tilde{p}\longrightarrow p$. First assume that $a_{i(p)}<s<a_{{i+1}(p)}$, this is, $\gamma_{i+1}(p)\subset{\Sigma}$ has internal tangencies

\begin{enumerate}
	\item [(a)] if $Z_{{i+1}(p)}=Z^s$ so $\gamma_{i+1}(p)\subset{\Sigma}$ and $q=\phi_{Z^s}(p,s)\in{\mathrm{A3}}$ with $q$ cusp visible point for $X$, Figure \ref{figure12} when $p_{i+1}\in{\Sigma^s}$. As $V$ a neighborhood of $M$ contain $p$ by inductive hypothesis $\tilde{s}\in{I_{Z_{{j+1}(\tilde{p})}}^{j+1}(\tilde{p})}$ with $Z_{{j+1}(\tilde{p})}=X$ for each $(\tilde{p},\tilde{s})\in{(V\times I)\cap \Omega_X}$ so
	$\displaystyle\lim_{\tilde{p}_j\to p_{i+1}} t_{Z^s}^+(\tilde{p}_j)=s-a_{i(p)}$, $\displaystyle\lim_{\tilde{p}\to p} \tilde{p}_{j+1}= \lim\limits_{\tilde{p}_j\to p_{i+1}} \varphi_{Z^s}(\tilde{p}_j,t_{Z^s}^+(\tilde{p}_j))=\varphi_{Z^s}(p_{i+1},s-a_{i(p)})=\phi_{Z^s}(p,s)=q$ and $\displaystyle\lim_{\tilde{p}\to p}a_{j(\tilde{p})}=\displaystyle\lim_{\tilde{p}\to p}a_{{j-1}(\tilde{p})} +\displaystyle\lim_{\tilde{p}_{j-1}\to p_{i+1}} t_{Z^s}^+(\tilde{p}_j)=a_{i(p)}+(s-a_{i(p)})=s$ and therefore
	
	and so
	\begin{align*}
	\lim\limits_{\substack{\tilde{p} \to p \\ \tilde{s}\to s}}  \phi_X(\tilde{p},\tilde{s}) &=\lim\limits_{\substack{\tilde{p} \to p \\ \tilde{s}\to s}}  \varphi_X(\tilde{p}_{j+1},\tilde{s}-a_{j(\tilde{p})})
	\\ &= \varphi_X(\lim_{\tilde{p}\to p} \tilde{p}_{j+1},\lim_{\tilde{s}\to s}\tilde{s}-\lim_{\tilde{p}\to p} a_{j(\tilde{p})})
	\\ &= \varphi_X(q,0)=q.
	\end{align*}

\begin{figure}[H]
	\includegraphics[scale=1.6]{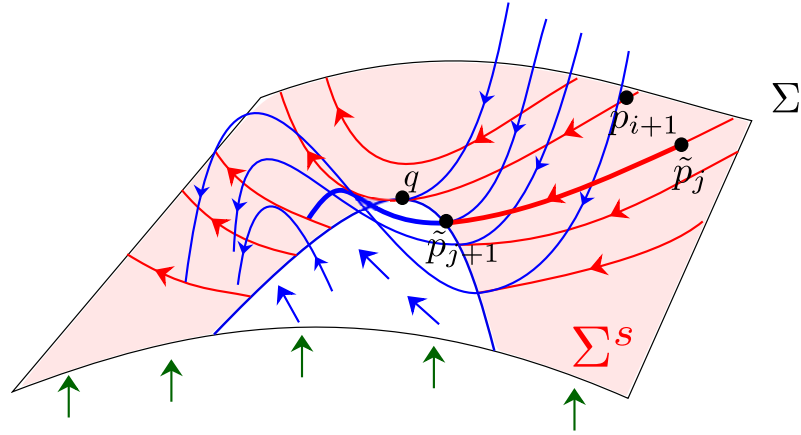}
	\caption{$q=\phi_Y(p,s)$ with $a_{i(p)}<s<a_{{i+1}(p)}$ is a point of case $\mathrm{A3}$, $\tilde{p}_j\in{\Sigma^{s}}$ and $\tilde{p}_{j+1}\in{\partial\Sigma^{s}}$.}
	\label{figure12}
\end{figure}
	\item [(b)] if $Z_{{i+1}(p)}=Y$ or $Z_{{i+1}(p)}=X$ so the process is analogous to item $(1)$ of the induction basis using the inductive hypothesis as in the below item (b).
\end{enumerate}
For $(p,s)\in{\overline{\Omega_X}}$ with $s=a_{{i+1}(p)}$ the process is analogous to item $(2)$ of the induction basis using the inductive hypothesis as in the item (b). Analogously for $(p,s)$ in $\overline{\Omega_Y}$ and $\overline{\Omega_{Z^s}}$.
\end{proof}
It is essential to highlight the following deduction.
\begin{cor}\label{coromegax}
The map $\phi_X:\overline{\Omega_X}\longrightarrow \Sigma^+$ is continuous and if $(p,s)\in{\overline{\Omega_X}\setminus \Omega_X }$ with $s\in{I_{Z_{i(p)}}^i(p)}$ then $\phi_X(p,s)=\phi_{Z_{i(p)}}(p,s)$. Analogously for the maps  $\phi_Y:\overline{\Omega_Y}\longrightarrow \Sigma^-$ and $\phi_{Z^s}:\overline{\Omega_{Z^s}}\longrightarrow \Sigma^S$.
\end{cor}

Now, we are ready to construct the semi-dynamical system for a $\mathsf{PSVF}$.
\begin{thm}\label{semiflow1}
	Assume $M$ is closed  3-dimensional $C^1$ manifold and $Z\in{\mathfrak{X}(M,h)}$. If $\Sigma=\Sigma^{c}\cup\Sigma^{s}\cup S_Z$ with $p \in{S_Z}$ is in some of the cases $\mathrm{A1}$ to $\mathrm{B2}$ then the trajectories of $Z=(X,Y)\in{\mathfrak{X}(M,h)}$ generate a semi-dynamical system $(M,\phi_{Z})$.
\end{thm}
For the proof of Theorem \ref{semiflow1}, we firstly show the following lemmas. 
\begin{lem}\label{sf1}
In the conditions of Proposition \ref{semiflow1}.

\begin{enumerate}
	\item [(a)] For all $(p,t)\in \overline{\Omega_X} \cap \overline{\Omega_{Z^s}}$ we have that $\phi_X(p,t)=\phi_{Z^s}(p,t)$;
	\item [(b)] For all $(p,t)\in \overline{\Omega_Y} \cap \overline{\Omega_{Z^s}}$  we have that $\phi_Y(p,t)=\phi_{Z^s}(p,t)$ and 
	\item [(c)] For all $(p,t)\in \overline{\Omega_X} \cap \overline{\Omega_Y}$  we have that $\phi_X(p,t)=\phi_Y(p,t)$.
\end{enumerate}
\end{lem}

\begin{proof}  Let $(p,t)\in \overline{\Omega_X} \cap \overline{\Omega_{Z^s}}$ then 
	
	\begin{itemize}
		\item [(1)] If $(p,t)\in{\Omega_X\cap\Omega_{Z^s}}$ then there exists $i\in{\Delta_Z^+(p)}$ such that $t\in{I_X^i(p)\cap I_{Z^s}^{i+1}(p)}$ (similarly for $t\in{I_{Z^s}^i(p)\cap I_{X}^{i+1}(p)}$) so $t=a_{i(p)}$ and		
	 $$\phi_X(p,t)=\varphi_X(p_{i+1},t-a_{i(p)})=\varphi_X(p_{i+1},0)=p_{i+1}=\varphi_{Z^s}(p_{i+1},0)=\phi_{Z^s}(p,t).$$
		
		\item [(2)] If $(p,t)\in (\overline{\Omega_X} \cap \overline{\Omega_{Z^s}})\setminus \Omega_X $ ( similarly for $(p,t)\in (\overline{\Omega_X} \cap \overline{\Omega_{Z^s}})\setminus \Omega_{Z^s} $) then there exists $i\in{\Delta_{Z^s}(p)}$ such that $t\in{I_{Z^s}^i(p)}$ and by Corollary \ref{coromegax} $\phi_X(p,t)=\phi_{Z^s}(p,t)$.
		
		\item [(3)] If $(p,t)\in (\overline{\Omega_X} \cap \overline{\Omega_{Z^s}})\setminus ( \Omega_X \cap \Omega_{Z^s})$ then $t=0$ and $p$ is a point os case $\mathrm{A1}$ so $\phi_X(p,t)=\varphi_X(p,0)=p=\varphi_{Z^s}(p,0)=\phi_{Z^s}(p,t)$.		
	\end{itemize}

The proof is analogous to $b)$ and $c)$.

\end{proof}

\begin{lem}\label{sf2}
With the same hypothesis that Proposition \ref{semiflow1}. The map  $\phi_Z:M \times \mathbb{R^+} \longrightarrow M$ such that 
	\begin{equation*}
	\phi_Z(p,t) = \left\{ \begin{array}{rcl}
	\phi_X(p,t)&& \text{if } (p,t)\in{\overline{\Omega_X}}, \\
	\phi_{Z^s}(p,t)&& \text{if } (p,t)\in{\overline{\Omega_{Z^s}}}, \\
	\phi_Y(p,t)&&\text{if }(p,t)\in{\overline{\Omega_Y}}
	\end{array}\right.
	\end{equation*}
	is well defined for all $(p,t)\in{M \times\mathbb{R^+}}$.
\end{lem}

\begin{proof}
For all $p\in{M}$, applying Definitions \ref{traj} and \ref{traj1} we have that
\begin{itemize}
	\item [(a)] if  $p\in{\textrm{int}(\Sigma^+)}$ then $\phi_Z(p,t)=\phi_X(p,t)$ for all $t\in{[0,a_{0(p)}]}$;
	\item [(b)] if $p\in{\textrm{int}(\Sigma^-)}$  then $\phi_Z(p,t)=\phi_X(p,t)$ for all $t\in{[0,a_{0(p)}]}$;
	\item [(c)] if  $p\in{\Sigma^{c+}}$ then $\phi_Z(p,t)=\phi_X(p,t)$ for all $t\in{[0,a_{0(p)}]}$;
	\item [(d)] if  $p\in{\Sigma^{c-}}$ then $\phi_Z(p,t)=\phi_Y(p,t)$ for all $t\in{[0,a_{0(p)}]}$;
	\item [(e)] if  $p\in{\Sigma^s}$ then $\phi_Z(p,t)=\phi_{Z^s}(p,t)$ for all $t\in{[0,a_{0(p)}]}$;
	\item [(f)] if $p$ is in the cases $\mathrm{A1}$, $\mathrm{A4}$ or $\mathrm{B1}$ that is fold or cusp visible point for $X$ then $\phi_Z(p,t)=\varphi_X(p,t)$ for all $t\in{[0,a_{0(p)}]}$,
	\item [(g)] if $p$ is in the case  $\mathrm{A2}$, $\mathrm{A3}$ or $\mathrm{B2}$ then $\phi_Z(p,t)=\phi_{Z^s}(p,t)$ for all $t\in{[0,a_{0(p)}]}$,
	\item [(h)] if $p$ is in the cases $\mathrm{A1}$, $\mathrm{A4}$ or $\mathrm{B1}$ that is fold or cusp visible point for $Y$ then $\phi_Z(p,t)=\varphi_Y(p,t)$ for all $t\in{[0,a_{0(p)}]}$.
\end{itemize}	
If $a_{0(p)}$ is a real number then $\phi_{Z_1(p)}(p,a_{0(p)})\in{\Sigma}$  hence apply item c) to h) above for all $t\in[a_{0(p)}, a_{1(p)}]$ and do it again if $a_1(p)=t_{Z_{0(p)}}^+(p_0)+t_{Z_{1(p)}}^+(p_1)<\infty$ and so forth. 	
\end{proof}

\begin{proof}[Proof of Theorem \ref{semiflow1}]
	
	Applying the gluing lemma for $\phi_X:\overline{\Omega_X}\longrightarrow M$, $\phi_Y:\overline{\Omega_Y}\longrightarrow M$ and $\phi_{Z^s}:\overline{\Omega_{Z^s}}\longrightarrow M$, by Remark \ref{rem3.14}, Proposition \ref{prop3.15}, Lemmas \ref{sf1} and \ref{sf2} we have that  $\phi_Z:M \times \mathbb{R^+} \longrightarrow M$ such that 
	\begin{equation*}
	\phi_Z(p,t) = \left\{ \begin{array}{rcl}
	\phi_X(p,t)&& \text{if } (p,t)\in{\overline{\Omega_X}}, \\
	\phi_{Z^s}(p,t)&& \text{if } (p,t)\in{\overline{\Omega_{Z^s}}}, \\
	\phi_Y(p,t)&&\text{if }(p,t)\in{\overline{\Omega_Y}}
	\end{array}\right.
	\end{equation*}
	is continuous.\\
	
	Now, we prove the items 1) and 2) of Definition \ref{semiflowdef}. The initial value property is satisfied since $\varphi_X$, $\varphi_{Z^s}$ and $\varphi_Y$ are flows. In order to prove the semi-group property, let $p\in{M}$ and $t,s\in{\mathbb{R^+}}$ and we prove that $\phi_Z(\phi_Z(p,t),s)=\phi_Z(p,t+s)$. Assume that $t\in{\I_{Z_{j(p)}}^j(p)}$ and let 	
	\begin{equation}
	q=\phi_Z(p,t)=\phi_{Z_{j(p)}}(p,t)=\varphi_{Z_{j(p)}}(p_j,t-a_{{j-1}(p)}),
	\end{equation}

	assume also that $s\in{I_{Z_{i(q)}}^i(q)}$ then  $t+s\in{I_{Z_{{j+i}(p)}}^{j+i}(p)}$ and $Z_{i(q)}=Z_{{j+i}(p)}$ so we use induction on $i\in{\Delta_Z^+(q)}$ for prove that 
	
	\begin{enumerate}
		\item $q_i=p_{j+i}$ for $i>0$,
		\item $a_{{i-1}(q)}=a_{{j+i-1}(p)}-t$ for $i>0$ and finally that
		\item $\phi_Z(\phi_Z(p,t),s)=\phi_Z(p,t+s)$.
	\end{enumerate}
	
	\textit{Induction basis:} Assume that $i=0$ then $s\in{I_{Z_{0(q)}}^0(q)}$,  $t+s\in{I_{Z_{j(q)}}^j(q)}$ and  $Z_{0(q)}=Z_{j(p)}$ and so  

	\begin{align*}
		\phi_Z(\phi_Z(p,t),s) &=\phi_Z(q,s)=\phi_{Z_{0(q)}}(q,s)
		\\ &= \varphi_{Z_{0(q)}}(\varphi_{Z_{j(p)}}(p_j,t-a_{j(p)}),s)
		\\ &= \varphi_{Z_{j(p)}}(\varphi_{Z_{j(p)}}(p_j,t-a_{j(p)}),s)
		\\ &= \varphi_{Z_{j(p)}}(p_j,t-a_{j(p)}+s)
		\\ &= \varphi_{Z_{j(p)}}(p_j,(t+s)-a_{j(p)})
		\\ &= \varphi_{Z_{j(p)}}(p,t+s)
		\\ &= \phi_Z(p,t+s).
	\end{align*}

	\textit{Inductive step:} For $i\in{\Delta_Z(q)}$ and $s\in{I_{Z_{i(q)}}^i(q)}$ assume true that  $q_i=p_{j+i}$, $a_{i-1}(q)=a_{j+i-1}(p)-t$ and $\phi_Z(\phi_Z(p,t),s)=\phi_Z(p,t+s)$ for all $t\in{\I_{Z_{j(p)}}^j(p)}$ and $s\in{\I_{Z_{i(q)}}^i(q)}$ then for $i+1\in{\Delta_Z(q)}$ and $s\in{I_{Z_{{i+1}(q)}}^{i+1}(q)}$
	
	\begin{enumerate}
		\item $q_{i+1}=\varphi_{Z_{i(q)}}(q_i,t_{Z_{i(q)}}^+(q_i))=\varphi_{Z_{{j+i}(p)}}(p_{j+i},t_{Z_{{j+i}(p)}}^+(p_{j+i}))= p_{j+i+1}$ ,
		\item $a_{i(q)}=a_{{i-1}(q)} + t_{Z_{i(q)}}^+(q_i)=(a_{{j+i-1}(p)}-t)+t_{Z_{{j+i}(p)}}^+(p_{j+i})=a_{{j+i}(p)}-t $ and so 
		\item  
		
		\begin{align*}
			\phi_Z(\phi_Z(p,t),s) &=\phi_Z(q,s)=\phi_{Z_{{i+1}(q)}}(q_{i+1},s-a_{i(q)})
			\\ &= \varphi_{Z_{{j+i+1}(p)}}(p_{j+i+1},s-(a_{{j+i}(p)}-t))
			\\ &= \varphi_{Z_{{j+i+1}(p)}}(p_{j+i+1},(t+s)-a_{{j+i}(p)})
			\\ &= \phi_{Z_{{j+i+1}(p)}}(p,t+s)
			\\ &= \phi_Z(p,t+s).
		\end{align*}
		
	\end{enumerate}
	
	Therefore $(M,\phi_Z)$ with $\Sigma=\Sigma^{c+}\cup \Sigma^{c-}\cup\Sigma^s\cup S_Z $ with $p \in{S_Z}$ is in some of the cases $\mathrm{A1}$ to $\mathrm{B2}$ is a semi-dynamical system. 
\end{proof}

As an outcome of Theorem \ref{semiflow1}, we obtain two significant corollaries.
\begin{cor}
	Assume $M$ is closed  2-dimensional $C^1$ manifold and $Z\in{\mathfrak{X}(M,h)}$. If $\Sigma=\Sigma^{c}\cup\Sigma^{s}\cup S_Z$ with $p \in{S_Z}$ is in some of the cases $\mathrm{A1}$, $\mathrm{A2}$ or $\mathrm{B1}$ then the trajectories of $Z=(X,Y)\in{\mathfrak{X}(M,h)}$ generate a semi-dynamical system $(M,\phi_{Z})$.
\end{cor}
	
\begin{figure}[H]
	\includegraphics[scale=1.4]{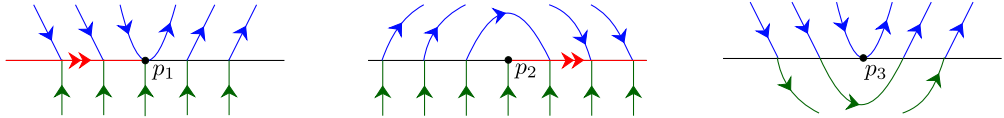}
	\caption{Example points of case $\mathrm{A}$ and $\mathrm{B}$ in 2-dimension.}
\end{figure}

\begin{cor}
		Assume $M$ is closed  3-dimensional $C^1$ manifold and $Z\in{\mathfrak{X}(M,h)}$. If $\Sigma=\Sigma^{c}\cup \mathrm{B1}$ then the trajectories of $Z=(X,Y)\in{\mathfrak{X}(M,h)}$ generate a dynamical system $(M,\varphi_{Z})$.
\end{cor}

\begin{proof}
Is sufficient note that in the proof of Theorem \ref{semiflow1} if no exists $\Sigma^s$ then we can work with $\varphi_X$ and $\varphi_Y$ instead of $\phi_X$ and $\phi_Y$ since, in this case, is possible to extend Definitions \ref{tx} and \ref{omegax} for $\mathbb{R^-}$. For $X \in{\mathfrak{X}(M)}$, let $\Lambda_X^{-}=\{p\in{\Sigma^+}; \ \varphi_X(p,(-\infty,0])  \nsubseteq  \Sigma^+\}$ and $t_X^+:\Sigma^+\longrightarrow \mathbb{R}^-\cup\{\infty\}$ such that
\begin{equation*}
t_X^-(p) = \left\{ \begin{array}{lll}
\ -\infty     &&  \text{if }  p\notin{\Lambda_X^-},\\	
\begin{array}{r@{}}
\mathrm{sup} \{t<0; \ \varphi_X(p,t)\in{S_X^{ic} \cup \Sigma_X^-} \}
\end{array}		
&&  \text{if } \varphi_X(p,[t,0])\subset{\Sigma^+} \text{ for } t\in{I\cap\mathbb{R}^-} \text{ with } t\neq 0,\\
\ 0  &&  \text{otherwise }  .
\end{array}\right.
\end{equation*}
Analogously for $Y\in{\mathfrak{X}(M)}$. If $p\in{M}$ then denoted $p_0=p$ and $p_i=		\varphi_{Z_{({i-1})(p)}^-}(p_{i-1},t_{Z_{{i-1}(p)}}^-(p_{i-1}))$ for $i\in{\mathbb{Z}^-}$, where  $Z_{0(p)}=X$ if $p\in{\mathrm{int}(\Sigma^+)}$, $Z_{0(p)}=Y$ if $p\in{\mathrm{int}(\Sigma^-)}$ and
\begin{equation*}
Z_{i(p)} = \left\{ \begin{array}{rcl}
X && \text{if }  p_{i}\in\Sigma^{c+} \text{or }  p_{i}\in{ \mathrm{B1}\cap S_X^v} ,\\
Y &&  \text{if } p_i\in{\Sigma^{c-}} \text{or } p_{i}\in{ \mathrm{B1}\cap S_Y^v}.
\end{array}\right.
\end{equation*}
Let $\Delta_X^-(p)=\{i\in{\mathbb{Z}^-}; \ Z_{i(p)}=X\}$ and $\Omega_X=\{(p,t)\in{M\times \mathbb{R}}; \ t\in{I_X(p)}\}$ whit $$I_X(p)=\bigcup_{i\in{\Delta_X^+ \cup \Delta_X^-}(p)} I_X^i(p)$$ 
where $I_X^0(p) = \left[ a_{0(p)}^-, a_{0(p)} \right] $,  $I_X^i(p) =	\left[ a_{{i-1}(p)}^-,a_{i(p)}^- \right] $ for $i<0$ such that $a_{i(p)}^-=\displaystyle\sum_{j=0}^{i} t_{Z_{j(p)}}^-(p_j)$. Analogously, we defined $\Omega_Y$ for the vector fields $Y$.

\end{proof}

Finally, the following theorem provides very general conditions under which an isolated invariant set must contain a periodic orbit in a $\mathsf{PSVF}$.

\begin{thm}\label{teo3.21}
	Assume $M$ is closed  3-dimensional $C^1$ manifold and $Z\in{\mathfrak{X}(M,h)}$. If $\Sigma=\Sigma^{c}\cup\Sigma^{s}\cup S_Z$ with $p \in{S_Z}$ is in some of the cases $\mathrm{A1}$ to $\mathrm{B2}$ and $\phi_Z:M\times [0,\infty)\rightarrow M$ is the semiflow generate by the trajectories of $Z=(X,Y)\in{\mathfrak{X}(M,h)}$. If $N$ is a isolating neighborhood for $\varphi$ which admits a Poincar\'{e} section $\varXi$ and either 
	\begin{equation}
	dim \ CH^{2n}(N,\varphi)=dim \ CH^{2n+1}(N,\varphi) \ \ \ \ \ for \ \ \ n\in{\mathbb{Z^+}}
	\end{equation}
	or 
	\begin{equation}
	dim \ CH^{2n}(N,\varphi)=dim \ CH^{2n-1}(N,\varphi) \ \ \ \ \ for \ \ \ n\in{\mathbb{Z^+}},
	\end{equation}
	where not all the above dimensions are zero, then $\varphi$ has a periodic trajectory in $N$.
\end{thm}

\begin{cor}\label{cor3.22}
	Under the hypotheses of Theorem \ref{teo3.21}, if $N$ has the Conley index of a hyperbolic periodic orbit, then $\mathrm{inv}(N)$ contains a periodic orbit.
\end{cor}


The following theorem is a n-dimensional version of Theorem 3.21, where we assume that the system has only crossing and sliding regions without tangency points.

\begin{thm}\label{teo3.23}
	Assume $M$ is closed  n-dimensional $C^1$ manifold and $Z\in{\mathfrak{X}(M,h)}$. If $\Sigma=\Sigma^{c}\cup\Sigma^{s}$ and $\phi_Z:M\times [0,\infty)\rightarrow M$ is the semiflow generate by the trajectories of $Z=(X,Y)\in{\mathfrak{X}(M,h)}$. If $N$ is a isolating neighborhood for $\varphi$ which admits a Poincar\'{e} section $\varXi$ and either 
	\begin{equation}
	dim \ CH^{2n}(N,\varphi)=dim \ CH^{2n+1}(N,\varphi) \ \ \ \ \ for \ \ \ n\in{\mathbb{Z^+}}
	\end{equation}
	or 
	\begin{equation}
	dim \ CH^{2n}(N,\varphi)=dim \ CH^{2n-1}(N,\varphi) \ \ \ \ \ for \ \ \ n\in{\mathbb{Z^+}},
	\end{equation}
	where not all the above dimensions are zero, then $\varphi$ has a periodic trajectory in $N$.
\end{thm}
\section{Some applications}

\subsection{Regularization}
 In this subsection, we provide an immediate consequence by the fact of the Conley index is robust under perturbation, we use Proposition \ref{prop2.11} and the concept of regularization of discontinuous vector fields, given by J. Sotomayor and  Marco A. Teixeira \cite{SotomayorTeixeira1}. By a transition function we mean a $C^{\infty}$ function $\varphi:\mathbb{R}\longrightarrow\mathbb{R}$ such that: $\varphi(t)=0$ if $t\leq -1$, $\varphi(t)=1$ if $t\geq 1$ and $\varphi^{\prime}>0$ if $t\in(-1,1)$.
\begin{definition}
A $\varphi_{\epsilon}-$regularization of $Z=(X,Y)$ is the one parameter family of vector fields $Z_{\epsilon}$ in $\mathfrak{X}(M)$ given by
$$Z_{\epsilon}(q)=\left(1- \varphi_{\epsilon}(h(q))\right) Y(q) + \varphi_{\epsilon}(h(q))X(q),$$
where $\varphi_{\epsilon}(t)=\varphi(\frac{t}{\epsilon})$.
\end{definition}

A fundamental characteristic of the Conley index is the homotopy invariance, and this permits the robustness of Theorem \ref{teo3.21}. Using the homotopy invariance of the Conley index and Proposition \ref{prop2.11} we have that the following proposition. 

\begin{prop}
Let $\gamma$ be a periodic orbit in a isolating neighborhood $N$ of the semiflow generated by the trajectories of $Z=(X,Y)$. Then there a $\epsilon_0$ such that for every $\epsilon\leq\epsilon_0$, $Z_{\epsilon}$ contains a periodic orbit in $N$. 
\end{prop}


\subsection{Closed poly-trajectories}
Piecewise-smooth vector fields defined on plane present a type of solutions called closed poly-trajectories, which are a particular case of periodic orbits defined in \cite{Sotomayor1}. This section provides some necessary conditions to guarantee the existence of closed poly-trajectories solutions in a disk of $M$ when the $\mathsf{PSVF}$ has sliding motion.

\begin{figure}[H]
	\includegraphics[scale=0.9]{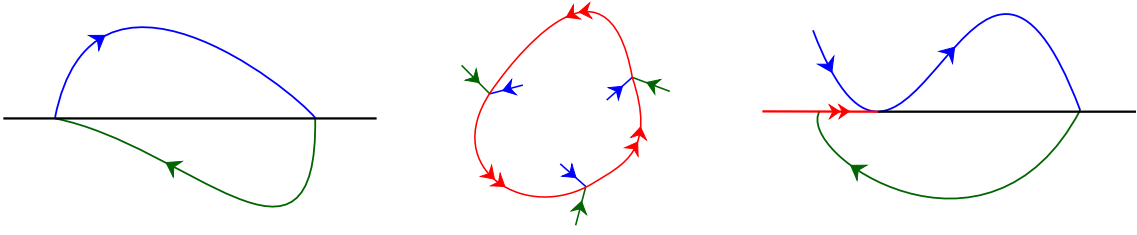}
	\caption{Closed poly-trajectories of kind I (left), II (center) and III (right).}
	\label{figure14}
\end{figure}

\begin{definition}
Consider $M$ is closed  2-dimensional $C^1$ manifold and $Z=(X,Y)\in{\mathfrak{X}(M,h)}$.
\begin{enumerate}
	\item A curve $\Gamma$ is a closed poly-trajectory if $\Gamma$ is closed and
	\begin{itemize}
		\item $\Gamma$ contains regular arcs of at least two of the vector fields $X|\Sigma^+$, $Y|\Sigma^-$, $Z^e$ and $Z^s$ or is composed by a single regular arc of either $Z^s$ or  $Z^e$;
		\item the transition between arcs of $X$ and arcs of $Y$ happens in sewing points (and vice versa);
		\item the transition between arcs of $X$ (or $Y$) and arcs of $Z^s$ or $Z^e$ happens through fold points or regular points in the escape or sliding arc, respecting the orientation. Moreover if $\Gamma\neq \Sigma$ then there exists at least one visible fold point on each connected component of $\Gamma\cap\Sigma$.
	\end{itemize}
\item Let $\Gamma$ be a canard cycle of $Z$. We say that	
	\begin{itemize}
		\item is a closed poly-trajectory of kind I if $\Gamma$ meets $\Sigma$ just in sewing
		points;
		\item is a closed poly-trajectory of kind II if $\Gamma=\Sigma$;
		\item is a closed poly-trajectory of kind III if $\Gamma$ contains at least one visible fold point of $Z$.
	\end{itemize}
	In Figures \ref{figure14}, \ref{figure15} and \ref{figure16} appear poly-trajectories of kind I, II and III respectively.
	\item Let $\Gamma$ be a closed poly-trajectory. We say that $\Gamma$ is hyperbolic if
	\begin{itemize}
		\item is of kind I and $\eta\prime(p)\neq 1$ where $\eta$ is the first return map defined on a segment $T$ with $p\in{T\cap\Gamma}$;
		\item is of kind II;
		\item $\Gamma$ is of kind III and either $\Gamma\cap\Sigma\subseteq\Sigma^c\cup\Sigma^s$ or $\Gamma\cap\Sigma\subseteq\Sigma^c\cup\Sigma^e$.
	\end{itemize}
	
\end{enumerate}
\end{definition}

Let $(pq)_{Z^*}$ be an arc of $Z^*$, where $Z^*=X,Y$, joining the visible fold point $p$ to the point $q\in{\Sigma_{Z^*}}$. We say that $(pq)_{Z^*}$ has \textit{focal kind} if there is not fold points between $p$ and $q$ (see Figure \ref{figure15} ) and we say that $(pq)_{Z^*}$ has \textit{graphic kind} if it has only one fold point between $p$ and $q$ (see Figure \ref{figure16}).

\begin{minipage}{\linewidth}
	\centering
	\begin{minipage}{0.35\linewidth}
		\begin{figure}[H]
			\includegraphics[width=\linewidth]{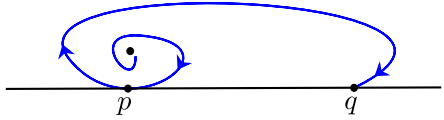}
			\caption{Focal kind arc.}
		    \label{figure15}
		\end{figure}
	\end{minipage}
	\hspace{0.05\linewidth}
	\begin{minipage}{0.35\linewidth}
		\begin{figure}[H]
			\includegraphics[width=\linewidth]{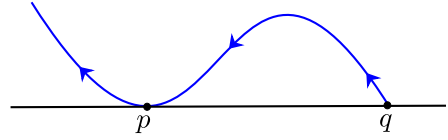}
			\caption{Graphic kind arc.}
			\label{figure16}
		\end{figure}
	\end{minipage}
\end{minipage}
So, we can use our main result to show the subsequent proposition. 

\begin{prop}\label{prop4.5}
Let $M$ a closed 2-dimensional $C^1$ manifold and $Z=(X,Y)\in{\mathfrak{X}(M,h)}$ such that $\Sigma=\Sigma^{c}\cup\Sigma^{s}\cup S_Z$ where for all $p \in{S_Z}$, $p$ is in some of the cases $\mathrm{A1}$, $\mathrm{A2}$ or $\mathrm{B1}$, moreover assume that $X$ is a linear vector field and  $Yh(x)>0$ for all $x\in \mathcal{U}\cap\Sigma$. Let $\mathcal{U}$ be a disk of $M$ that belong to the same chart of $M$ and such that:
\begin{itemize}
	\item [(1)] $\textrm{int}(\mathcal{U})\cap \Sigma\neq\emptyset$,
	\item [(2)]	$\mathcal{U}\cap\textrm{int}(\Sigma^+)$ contain only one equilibrium point $\tilde{x}$, is which unstable focus and an the unstable manifold of $\tilde{x}$ intercept to $\Sigma$ is an arc in $\Sigma$, and
	\item [(3)] there are no pseudo-equilibrium points in $\mathcal{U}\cap\Sigma$.
\end{itemize}	
Then, in $\mathcal{U}$ there exists a hyperbolic poly-trajectory of kind $III$.
\end{prop}

\begin{proof} A closed poly-trajectory of $Z=(X,Y)\in{\mathfrak{X}(M,h)}$ is a periodic orbit for semi-flow $\phi_Z$ thus, we used Theorem \ref{teo3.21}.\newline

Let $\tilde{x}$ the unstable focus in $\mathcal{U}\cap\textrm{int}(\Sigma^+)$, $p$ the visible fold point for $X$ such that $\gamma_0(p)=p\cdot[0,t_X^+(p)]$ is the focal kind arc. Let   $\tilde{p}$ and $\tilde{q}$ the intersection of $\partial \mathcal{U}$ and $\Sigma$ such that $Xh(x)>0$ for all $x\in{(\tilde{p}p)_{\Sigma}\setminus{\{p\}}}$ where $(\tilde{p}p)_{\Sigma}$ is the arc in $\Sigma$ of the point $\tilde{p}$ to $p$, and $Xh(x)<0$ for all $x\in{(p\tilde{q})_{\Sigma}\setminus{\{p\}}}$. Take $q\in{\textrm{int}(\tilde{p}p)_{\Sigma}}$ such that $q_1\in{\textrm{int}(p\tilde{q})_{\Sigma}}$ and $\gamma_0(q)\subset{\mathcal{U}}$. Consider $\mu <0$ such that $\Sigma_{\mu}=h^{-1}(\mu)$ is parallel to $\Sigma$ and $\Sigma_{\mu}\cap\textrm{int}(\mathcal{U})\neq \emptyset$, and consider the points $v, \tilde{v}\in{\Sigma_{\mu}}$ satisfying that $q=\tilde{v}_1$ and $q_1=v_1$. Let $\tilde{N}$ the region bounded by the curve $(v\tilde{v})_{\Sigma_{\mu}}\cup\gamma_0(\tilde{v})\cup\gamma_0(q) \cup \gamma_0(v)$, see Figure \ref{figure17}.
\begin{figure}[H]
	\includegraphics[scale=0.9]{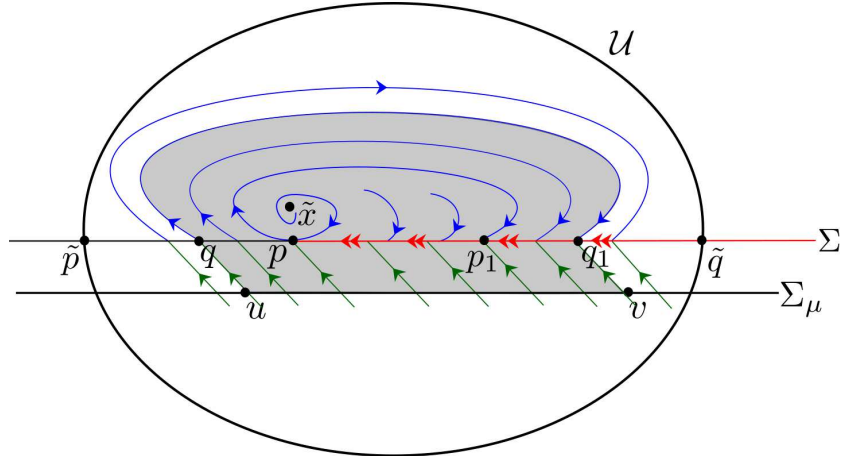}
	\caption{Construction of the isolating neighborhood.}
	\label{figure17}
\end{figure}
Now, let $\Xi$ a local section for $\varphi_X$ crossing the stable manifold of $\tilde{x}$ but $\varphi_X\cap\Sigma=\emptyset$, as in the Figure \ref{figure18}. If $a\pm bi$ are the eigenvalues of $X$ associate to $\tilde{x}$ then, let $r\in{\Xi}$ such that $\theta(r)<90^\circ$ when $b<0$ or $90^\circ<\theta(r)<180^\circ$ when $b>0$, where $\theta(r)$ is the angle that the vector $X(r)$ makes with $\Sigma$. As $\tilde{x}$ is a unstable focus, then there exist $\tilde{r}\in{\Xi}$ and let $\tilde{\gamma}_{r}$ the arc of $r$ to $\tilde{r}$ by flow $\varphi_X$, then $\tilde{\gamma}_{r}\cup (r\tilde{r})_{\Xi}$ is homeomorphic to $S^1$. Thus, let $N=\tilde{N}\setminus \left(\tilde{\gamma}_{r}\cup (r\tilde{r})_{\Xi}\right) $.
\begin{figure}[H]
	\includegraphics[scale=1.2]{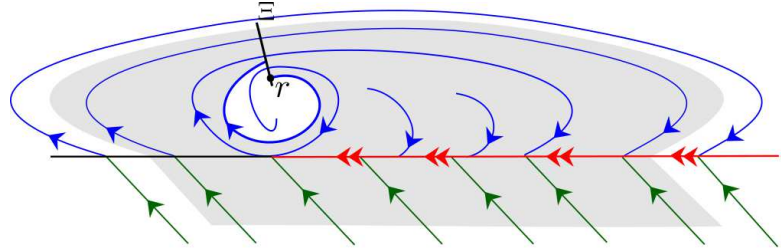}
	\caption{Poincaré Section.}
	\label{figure18}
\end{figure}
Note that $\mathrm{inv}(N)=\gamma_0(p)\cup (pp_1)_{\Sigma}$ and so $\mathrm{inv}(N)\subset \mathrm{int}(N)$. The exit set $L$ is empty, in fact, if $x\in N\cap \Sigma^+$ then by hypotheses (3) either $x\in{\Sigma^s}$ or there exists $t_X^+(x)\in{(0,\infty)}$ such that $x \cdot t_X^+(x)\in{\Sigma^s}$, and by hypotheses (5), there is $s_x\geq 0$ such that $x \cdot (t_X^+(x) +s_x)\in{\mathrm{inv}(N)}$; analogously if $x\in N\cap \Sigma^-$. In the Figure \ref{figure18} you can see that the homotopy type of $N$ is the same of a stable periodic orbit thus
\begin{equation*}
CH^k(N) \approx \left \{ \begin{matrix} \mathbb{Z} & \mbox{ }k=0, 1,
\\ 0 &  \mbox{ otherwise}.
\end{matrix}\right. 
\end{equation*}

\begin{figure}[H]
	\includegraphics[scale=0.6]{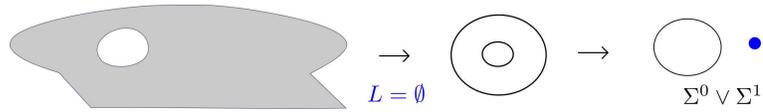}
	\caption{The homotopy type of $N/L$.}
	\label{figure19}
\end{figure}
    \textit{Poincaré Section.} Let $\Xi$ a local section for $\varphi_X$ crossing the stable manifold of $\tilde{x}$ as in the Figure \ref{figure19}. We claim that $\Xi$ is the required Poincare section for $N$. It is closed, and it is transverse to the semiflow. Finally, we must show that the forward orbit of every point in $N$ intersects $\Xi$.  If $x\in{N\cap\Sigma^+}$ then, by the hypothesis (3), $x\in{\Lambda_X^{+}}$ and there is $t_X^+(x)>0$ such that $x_1= x\cdot t_X^+(x)\in{\Sigma\cap N}$, moreover by the hypothesis (a), $x_1\in{\Sigma^s}$. By hypothesis (5), $x_1\in{\Lambda_{Z^+}^{+}}$ then there is $t_{Z^s}^+(x_1)>0$ such that $x_2= x_1\cdot t_{Z^s}^+(x_1)=p$ and forward orbit of $p$ intersects $\Xi$. If $x\in{N\cap\Sigma^-}$ then,  by the hypothesis (4) and by construction of $N$, there is $t_Y^+(x)>0$ such that $x_1= x\cdot t_Y^+(x)\in{\Sigma\subset\Sigma^+}$ and continue as before. Thus $\Xi$ is a Poincare section for the semiflow $\phi_Z$ in $N$.       
%
\end{proof}
Using the ideas of Proposition \ref{prop4.5}, we can show that in a disk $\mathcal{U}$ with dynamic as in the figure \ref{figure20}, there exists a hyperbolic poly-trajectory of kind $III$.

\begin{figure}[H]
	\includegraphics[scale=0.85]{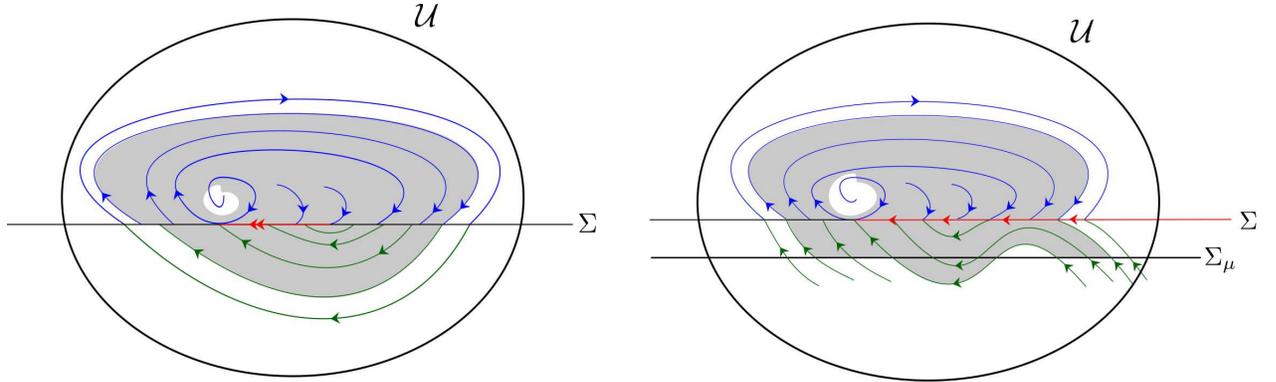}
	\caption{In both cases, there exists a hyperbolic poly-trajectory of kind $III$ contain to $\mathcal{U}$.}
	\label{figure20}
\end{figure}

\section{Closing remarks and future directions}

The results in this paper bring us one step closer to generalizing the tool of  C. Mccord, K. Mischaikow, and M. Mrozek in \cite{Mccord} to piecewise smooth vector fields using the Filippov convention. The next step is to include tangential singularities with orders greater than $3$ in systems defined in a manifold of dimension $n>3$. Forthcoming works include applying the results to guarantee the existence of periodic orbits in a piecewise smooth system that models an intermittent treatment of the human immunodeficiency virus, see \cite{de2020global}. Another application is to find periodic orbits in hill models in biology using $\mathsf{PSVF}$ to defined switch systems as in Dynamic Signatures Generated by Regulatory Networks ($\mathsf{DSGRN}$), see \cite{cummins2016combinatorial}.

\section*{Acknowledgment}

The first author is is supported by CAPES (the Coordena{\c c}{\~a}o Aperfei{\c c}oamento de Pessoal de N{\'i}vel Superior-Brasil) and affiliated with IME-UFG (Instituto de Matem\'atica e Estat\'istica, Universidade Federal de Goi\'as). The second author  affiliated with DIMACS (the Center for Discrete Mathematics and Theoretical Computer Science), Rutgers University,  and IME-UFG (Instituto de Matem\'atica e Estat\'istica, Universidade Federal de Goi\'as) and would like to acknowledge the support of the National Science Foundation under grant HDR TRIPODS 1934924.

\nocite{*}
\bibliographystyle{abbrv}  

\bibliography{references}

\begin{thebibliography}{10}

\bibitem{bernardo}
M.~Bernardo, C.~Budd, A.~R. Champneys, and P.~Kowalczyk.
\newblock {\em Piecewise-smooth dynamical systems: theory and applications},
  volume 163.
\newblock Springer Science \& Business Media, 2008.

\bibitem{carmona}
V.~Carmona, E.~Freire, E.~Ponce, and F.~Torres.
\newblock Bifurcation of invariant cones in piecewise linear homogeneous
  systems.
\newblock {\em International Journal of Bifurcation and Chaos},
  15(08):2469--2484, 2005.

\bibitem{Casagrande}
R.~Casagrande, K.~de~Rezende, and M.~Teixeira.
\newblock The conley index for discontinuous vector fields.
\newblock {\em Geometriae Dedicata}, 136(1):47, 2008.

\bibitem{Conley}
C.~C. Conley.
\newblock {\em Isolated invariant sets and the Morse index}.
\newblock Number~38. American Mathematical Soc., 1978.

\bibitem{cummins2016combinatorial}
B.~Cummins, T.~Gedeon, S.~Harker, K.~Mischaikow, and K.~Mok.
\newblock Combinatorial representation of parameter space for switching
  networks.
\newblock {\em SIAM journal on applied dynamical systems}, 15(4):2176--2212,
  2016.

\bibitem{de2020global}
T.~de~Carvalho, R.~Cristiano, L.~F. Gon{\c{c}}alves, and D.~J. Tonon.
\newblock Global analysis of the dynamics of a mathematical model to
  intermittent hiv treatment.
\newblock {\em Nonlinear Dynamics}, 101(1):719--739, 2020.

\bibitem{doedel2003computation}
E.~J. Doedel, R.~C. Paffenroth, H.~B. Keller, D.~Dichmann, J.~Gal{\'a}n-Vioque,
  and A.~Vanderbauwhede.
\newblock Computation of periodic solutions of conservative systems with
  application to the 3-body problem.
\newblock {\em International Journal of Bifurcation and Chaos},
  13(06):1353--1381, 2003.

\bibitem{du}
Z.~Du, Y.~Li, and W.~Zhang.
\newblock Bifurcation of periodic orbits in a class of planar filippov systems.
\newblock {\em Nonlinear Analysis: Theory, Methods \& Applications},
  69(10):3610--3628, 2008.

\bibitem{euzebio2014estudo}
R.~D. Euz{\'e}bio.
\newblock Estudo de conjuntos minimais para sistemas descont{\'\i}nuos em
  dimens{\~o}es 2 e 3.
\newblock 2014.

\bibitem{Filippov}
A.~Filippov.
\newblock {\em Differential equations with discontinuous right-hand sides}.
\newblock Mathematics and its Applications (Soviet Series), Kluwer Academic
  Publishers-Dordrecht, 1988.

\bibitem{Teixeira}
O.~M. Gomide and M.~A. Teixeira.
\newblock On structural stability of 3d filippov systems.
\newblock {\em Mathematische Zeitschrift}, 294(1-2):419--449, 2020.

\bibitem{junior2016orbitas}
C.~J{\'u}nior and R.~Pazim.
\newblock {\'O}rbitas peri{\'o}dicas em sistemas diferenciais suaves por
  partes.
\newblock 2016.

\bibitem{li}
T.~Li and X.~Chen.
\newblock Periodic orbits of linear filippov systems with a line of
  discontinuity.
\newblock {\em Qualitative Theory of Dynamical Systems}, 19(1):1--22, 2020.

\bibitem{llibre2}
J.~Llibre, A.~C. Mereu, and D.~D. Novaes.
\newblock Averaging theory for discontinuous piecewise differential systems.
\newblock {\em Journal of Differential Equations}, 258(11):4007--4032, 2015.

\bibitem{llibre3}
J.~Llibre and E.~Ponce.
\newblock Three nested limit cycles in discontinuous piecewise linear
  differential systems with two zones.
\newblock {\em Dyn. Contin. Discrete Impuls. Syst. Ser. B Appl. Algorithms},
  19(3):325--335, 2012.

\bibitem{llibre}
J.~Llibre and M.~A. Teixeira.
\newblock Periodic orbits of continuous and discontinuous piecewise linear
  differential systems via first integrals.
\newblock {\em S{\~a}o Paulo Journal of Mathematical Sciences}, 12(1):121--135,
  2018.

\bibitem{llibre1}
J.~Llibre, D.~J. Tonon, and M.~Q. Velter.
\newblock Crossing periodic orbits via first integrals.
\newblock {\em International Journal of Bifurcation and Chaos}, 30(11):2050163,
  2020.

\bibitem{Mccord}
C.~McCord, K.~Mischaikow, and M.~Mrozek.
\newblock Zeta functions, periodic trajectories, and the conley index.
\newblock {\em Journal of differential equations}, 121(2):258--292, 1995.

\bibitem{Mischaikow}
K.~Mischaikow.
\newblock Conley index theory.
\newblock In {\em Dynamical systems}, pages 119--207. Springer, 1995.

\bibitem{Mrozek}
M.~Mrozek.
\newblock The conley index on compact anr’s is of finite type.
\newblock {\em Results in Mathematics}, 18(3-4):306--313, 1990.

\bibitem{mrozek2020creating}
M.~Mrozek and T.~Wanner.
\newblock Creating semiflows on simplicial complexes from combinatorial vector
  fields.
\newblock {\em arXiv preprint arXiv:2005.11647}, 2020.

\bibitem{Rybakowski}
K.~P. Rybakowski.
\newblock {\em The homotopy index and partial differential equations}.
\newblock Springer Science \& Business Media, 1987.

\bibitem{Abel}
A.~Sagodi.
\newblock Conley index theory in neuroscience.
\newblock Master Thesis, University of Amsterdam, 2020.

\bibitem{Salamon}
D.~Salamon.
\newblock Connected simple systems and the conley index of isolated invariant
  sets.
\newblock {\em Transactions of the American Mathematical Society},
  291(1):1--41, 1985.

\bibitem{Sotomayor1}
J.~Sotomayor and A.~L.~F. Machado.
\newblock Structurally stable discontinuous vector fields in the plane.
\newblock {\em Qualitative Theory of Dynamical Systems}, 3(1):227--250, 2002.

\bibitem{SotomayorTeixeira1}
J.~Sotomayor and M.~Teixeira.
\newblock Regularization of discontinuous vector fields.
\newblock In {\em Proceedings of the international conference on differential
  equations, Lisboa}, pages 207--223. World Scientific, 1996.

\bibitem{SotomayorTeixeira}
J.~Sotomayor and M.~A. Teixeira.
\newblock Vector fields near the boundary of a 3-manifold.
\newblock In {\em Dynamical Systems Valparaiso 1986}, pages 169--195. Springer,
  1988.

\bibitem{Sotomayor}
T.~Sotomayor et~al.
\newblock Structural stability manifolds with boundary.
\newblock In {\em Global analysis and its applications}. 1974.

\bibitem{Teixeira1}
M.~A. Teixeira.
\newblock Stability conditions for discontinuous vector fields.
\newblock {\em Journal of Differential Equations}, 88(1):15--29, 1990.

\bibitem{Cameron2}
C.~Thieme.
\newblock Isolating neighborhoods and filippov systems: Extending conley theory
  to differential inclusions.
\newblock {\em arXiv preprint arXiv:1912.13116}, 2019.

\bibitem{Cameron1}
C.~Thieme.
\newblock Multiflows: A new technique for filippov systems and differential
  inclusions.
\newblock {\em arXiv preprint arXiv:1905.07051}, 2019.

\bibitem{Cameron3}
C.~Thieme.
\newblock Conley index theory and the attractor-repeller decomposition for
  differential inclusions.
\newblock {\em arXiv preprint arXiv:2009.00696}, 2020.

\bibitem{tonon2010sistemas}
D.~J. Tonon.
\newblock Sistemas de filippov em variedades tridimensionais.
\newblock 2010.

\end{thebibliography}

\end{document}